\documentclass[a4paper, final, 10pt]{article}

\usepackage{amsfonts,amsbsy,amssymb,amsthm}
\usepackage{amssymb}
\usepackage{amsthm}

\usepackage{makeidx}
\usepackage{mathtools}
\usepackage{bm}
\usepackage{cite}
\usepackage{showkeys}
\usepackage[english]{babel}
\usepackage{dblaccnt}
\usepackage{accents}
\usepackage{graphicx}
\usepackage{psfrag}
\usepackage{subfig}
\usepackage{color}
\usepackage{float}
\usepackage{amscd}
\usepackage[mathscr]{euscript}
\usepackage{lipsum}
\usepackage{epstopdf}
\usepackage{algorithm}
\usepackage{algpseudocode}

\newtheorem{remark}{Remark}[section]

\numberwithin{equation}{section}

\begin{document}

\title{Convexification method for a coefficient inverse problem and its
performance for experimental backscatter data for buried targets\thanks{Supported 
by US Army Research Laboratory and US Army Research Office grant
W911NF-15-1-0233 and by the Office of Naval Research grant N00014-15-1-2330.
In addition, the work of Kolesov A.E. was partially supported by Mega-grant
of the Russian Federation Government (N14.Y26.31.0013) and RFBR (N17-01-00689A)}}
\author{ Michael V. Klibanov\thanks{Department of Mathematics \& Statistics, University of North Carolina at
Charlotte, Charlotte, NC 28223, USA (mklibanv@uncc.edu)} \and Aleksandr E.
Kolesov\footnotemark[2] \thanks{Institute of Mathematics and Information Science, North-Eastern Federal
University, Yakutsk, Russia (akolesov@uncc.edu)} \and Dinh-Liem Nguyen%
\thanks{Department of Mathematics, Kansas State University, Manhattan, KS 66506
(dlnguyen@ksu.edu)} }
\date{}
\maketitle

\begin{abstract}
We present in this paper a novel numerical reconstruction method for solving
a 3D coefficient inverse problem with scattering data generated by a single
direction of the incident plane wave. This inverse problem is well-known to
be a highly nonlinear and ill-posed problem. Therefore, optimization-based
reconstruction methods for solving this problem would typically suffer from
the local-minima trapping and require strong a priori information of the
solution. To avoid these problems, in our numerical method, we aim to
construct a cost functional with a globally strictly convex property, whose
minimizer can provide a good approximation for the exact solution of the
inverse problem. The key ingredients for the construction of such functional
are an integro-differential formulation of the inverse problem and a
Carleman weight function. Under a (partial) finite difference approximation,
the global strict convexity is proven using the tool of Carleman estimates.
The global convergence of the gradient projection method to the exact
solution is proven as well. We demonstrate the efficiency of our
reconstruction method via a numerical study of experimental backscatter data
for buried objects.
\end{abstract}

\textbf{Keywords.} Carleman weight function, Carleman estimates, reconstruction method,
convexification, global convergence, coefficient inverse problem,
experimental  data

\textbf{AMS subject classification.} 35R30, 78A46, 65C20

\section{Introduction}

\label{sec:1} We develop in this paper a novel numerical method for solving
a coefficient inverse problem (CIP) for the 3D Helmholtz equation with
scattering data generated by a single direction of the incident plane wave
at multiple frequencies. More precisely, the goal of this CIP is to recover
a coefficient in the Helmholtz equation from boundary measurements of its
solutions for a single direction of the incident plane wave at multiple
frequencies.

This CIP arises in a wide range of applications including non-destructive
testing, detection of explosives, medical imaging, geophysics, etc. It is
also well-known that any CIP is a highly nonlinear and ill-posed problem
causing substantial challenges in the design of numerical algorithms for
solving it. Optimization-based reconstruction methods can be considered as
the most studied approach for solving CIPs in general. However, these
methods suffer from the fact that they might converge to a local minimum,
which is not the true solution of the CIP. Moreover, these methods typically
require strong a priori information of the solution, which is not always
available in practice.

The goal of the so-called \emph{globally convergent method} (GCM), recently
developed by the first author and coauthors (see e.g.~\cite%
{BeilinaKlibanov12}) is to overcome the drawbacks mentioned above when
solving CIPs. This method aims to provide a point in a sufficiently small
neighborhood of the true solution of the CIP without any advanced knowledge
of this neighborhood. The size of this neighborhood should depend only on
approximation errors and the level of noise in the data.

The numerical method we develop in this paper can be considered as the
second type of GCMs, which has certain advantages compared with the first
type of GCMs in \cite%
{BeilinaKlibanov12,KlibanovLiem17buried,KlibanovLiem17exp}. More precisely,
we do not impose in the convergence analysis here the assumption on a small
interval of wavenumbers. Neither we do not iterate here with respect to the
so-called tail functions. The combination of the latter two features with
the globally strictly convex cost functional (below) are the main
improvements of the convexification over the first type of globally
convergent methods.

This second type of GCMs is also called \emph{convexification} methods,
which was studied for the 1D case in~\cite{KlibanovKolesov17}. The present
work can be considered as a generalization to the 3D case of the cited 1D
version. Convexification methods are based on the minimization of the
weighted cost functional with a Carleman weight function (CWF) in it. The
CWF is the function which is involved in the Carleman estimate for the
corresponding PDE operator. The CWF can be chosen in such a way that the
cost functional becomes strictly convex. Note that the majority of known
numerical methods of solutions of nonlinear ill-posed problems minimize
conventional least squares cost functionals (see, e.g. \cite%
{Chavent09,Goncharsky13,Goncharsky17}), which are usually non convex and
have multiple local minima and ravines, see, e.g. \cite{Scales92} for a good
numerical example of multiple local minima.

We work in this paper with a semidiscrete version of the convexification,
which is more realistic for computations than continuous versions used in
previous works on the convexification of the first author with coauthors 
\cite%
{BakushinskiiKlibanov17,Klibanov97a,Klibanov15,KlibanovKolesov17,KlibanovThanh15,KlibKol3D}%
. \textquotedblleft Semidiscrete" means that we develop the theory for the
case when the differential operator we work with is written in finite
differences with respect to two out of three spatial variables and in the
continuous form with respect to the third variable. We impose a
computationally reasonable assumption that the grid step size does not tend
to zero (unlike the case of some forward problems). The fully discrete case,
i.e. when derivatives with all three variables are written in finite
differences, is not investigated yet. Indeed, it is well known that this
case is quite a complicated one for ill-posed problems for PDEs, especially
in nonlinear cases, such as we work with. There are known only a few results
for the fully discrete cases of linear ill-posed problems, see, e.g. \cite%
{Burman,KS}. We also refer to the recent publication \cite{KlibKol3D} of the
first two authors about a 3D version of the convexification method. In 
\cite{KlibKol3D} convexification was numerically tested on some
computationally simulated data. This is unlike the current paper in which
testing is done for a significantly more challenging case of experimental
data. The theory in \cite{KlibKol3D} is developed for the continuous case.
Although the idea of the semidiscrete version is briefly outlined in \cite%
{KlibKol3D}, corresponding theorems are neither formulated nor proved there,
unlike the current paper.

We point out that the CIP considered in this paper is also called a inverse
scattering problem in some contexts. There is a vast literature on both
theoretical and numerical studies on this inverse problem and its
variations, see, e.g. \cite%
{Am1,Am2,Cakon2006,Colto2013,It,Ito,Kab1,Kab2,Lakhal1,Liu1,Liu2,Kar1,Kar2}.
These cited papers have considered the cases of multiple measurements and/or
shape reconstructions. We recall that we consider in this paper the CIP with
a single measurement which is both different and more challenging than the
configurations considered in those cited papers.

In the next section, we provide a statement of the forward and inverse
problems. In Section 3 we present an integro-differential equation
formulation of the CIP. Section 4 involves the approximation of the tail
function which is an important component in the integro-differential
equation. We introduce in Section 5 the partial finite difference
approximation and related function spaces for the integro-differential
formulation. We describe in Section 6 the weighted cost functional with the
Carleman weight function in it. Section 7 is dedicated to the theoretical
analysis, including a Carleman estimate and proofs of global strict
convexity of that functional as well as convergence results for the
optimization problem. Finally, our numerical study is presented in Section 8.

\section{Problem Statement}

\label{sec:2.2}

Let $\mathbf{x}=(x,y,z)\in \mathbb{R}^{3}$ and consider positive numbers $%
b>0 $ and $d>0$. For the convenience for our numerical study (Section 8), we
define from the beginning the domain of interest $\Omega $ and the
backscatter part $\Gamma $ of its boundary as 
\begin{equation}
\begin{gathered} \Omega =\left\{ \left( x,y,z\right) :\left\vert
x\right\vert, \left\vert y\right\vert <b, z\in \left( -\xi ,d\right)
\right\} , \Gamma =\left\{ (x,y,z):|x|,|y|<b,z=-\xi \right\} .
\label{eq:2.1} \end{gathered}
\end{equation}%
Let the function $c(\mathbf{x})$ be the spatially distributed dielectric
constant and $k$ be the wavenumber. We consider the forward scattering
problem for the Helmholtz equation: 
\begin{align}  \label{eq:helmholtz}
&\Delta u+k^{2}c(\mathbf{x})u=0,\quad \mathbf{x}\in \mathbb{R}^{3}, \\
&u(\mathbf{x},k)=u_{s}(\mathbf{x},k)+u_{i}(\mathbf{x},k), \\
&\lim_{r\rightarrow \infty }r\left(\partial u_{s}/\partial r
-iku_{s}\right) =0,\quad r=\left\vert \mathbf{x}\right\vert,  \label{eq:2.3}
\end{align}
where $u(\mathbf{x},k)$ is the total wave, $u_{i}(\mathbf{x},k)$ is the
incident wave and $u_{s}(\mathbf{x},k)$ is the scattered wave satisfying the
Sommerfeld radiation condition. This condition means that the scattered
field behaves like a outgoing spherical wave far away from the scattering
medium.

Here we consider $u_{i}(\mathbf{x},k)$ as the incident plane wave
propagating along the positive direction of the $z-$axis: 
\begin{equation}
u_{i}(\mathbf{x},k)=e^{ikz}.  \label{eq:uinc}
\end{equation}%
%
%
%
%
Also, the function $c(\mathbf{x})$ satisfies with the following conditions: 
\begin{equation}
c(\mathbf{x})=1+\beta (\mathbf{x}),\quad \beta (\mathbf{x})\geq 0,\,\mathbf{x%
}\in \mathbb{R}^{3},\quad \mbox{and }c(\mathbf{x})=1,\,\mathbf{x}\notin 
\overline{\Omega }.  \label{eq:coef}
\end{equation}%
The assumption of (\ref{eq:coef}) $c(\mathbf{x})=1$ in $\mathbb{R}%
^{3}\setminus \Omega $ means that we have vacuum outside of the domain $%
\Omega .$ Finally, we assume that $c(\mathbf{x})\in C^{15}(\mathbb{R}^{3})$.
This smoothness condition was imposed to derive the asymptotic behavior of
the solution of the Helmholtz equation (\ref{eq:helmholtz}) (see \cite%
{KlibanovRomanov16}). We also note that extra smoothness conditions are
usually not of a significant concern when a CIP is considered, see, e.g.
Theorem 4.1 in \cite{Rom3}. 
Also, it follows from Lemma 3.3 of \cite{KlibanovLiem16} that the derivative 
$\partial _{k}u(\mathbf{x},k)$ exists for all $\mathbf{x}\in \mathbb{R}%
^{3},k>0$ and satisfies the same smoothness condition as the function $u(%
\mathbf{x},k).$

\textbf{Coefficient Inverse Problem (CIP).} \emph{Let $\Omega $ and $\Gamma
\subset \partial \Omega $ be as in (\ref{eq:2.1}). Let the wavenumber $k\in
\lbrack \underline{k},\overline{k}],$ where }$[\underline{k},\overline{k}]$%
\emph{$\subset \left( 0,\infty \right) $ is an interval of wavenumbers.
Determine the function $c(\mathbf{x}),\,\mathbf{x}\in \Omega $, given the
boundary data $g_{0}(\mathbf{x},k)$ as} 
\begin{equation}
u(\mathbf{x},k)=g_{0}(\mathbf{x},k),\quad \mathbf{x}\in \Gamma ,\,k\in
\lbrack \underline{k},\overline{k}].  \label{eq:cisp}
\end{equation}

In addition to the data (\ref{eq:cisp}) we can obtain the boundary
conditions for the derivative of the function $u(\mathbf{x},k)$ in the $z-$%
direction using the data propagation procedure (see~\cite%
{KlibanovLiem17buried}), 
\begin{equation}
u_{z}(\mathbf{x},k)=g_{1}(\mathbf{x},k),\quad \mathbf{x}\in \Gamma ,\,k\in
\lbrack \underline{k},\overline{k}].  \label{eq:gz0}
\end{equation}%
Even though we use the data propagation procedure in our computations below,
we do not describe it here for brevity. Instead, we refer to detailed
descriptions in \cite{KlibanovLiem17buried, KlibanovLiem17exp}. In fact,
this procedure is widely used in Optics under the name the \emph{angular
spectrum representation}.

In addition, we complement Dirichlet (\ref{eq:cisp}) and Neumann (\ref%
{eq:gz0}) boundary conditions on $\Gamma $ with the heuristic Dirichlet
boundary condition at the rest of the boundary $\partial \Omega $ as:%
\begin{equation}
u(\mathbf{x},k)=e^{ikz},\quad \mathbf{x}\in \partial \Omega \setminus \Gamma
,k\in \lbrack \underline{k},\overline{k}].  \label{2.2}
\end{equation}%
The boundary condition (\ref{2.2}) coincides with the one for the uniform
medium with $c\left( \mathbf{x}\right) \equiv 1.$ To justify (\ref%
{eq:helmholtz}), we recall that, using the tail functions method, it was
demonstrated in sections 7.6 and 7.7 of \cite{KlibanovLiem17exp} that (\ref%
{eq:helmholtz}) does not affect much the reconstruction accuracy as compared
with the correct Dirichlet boundary condition. Besides, (\ref{eq:helmholtz})
has always been used in works \cite{KlibanovLiem17exp} with experimental
data, where accurate results were obtained by the tail functions globally
convergent method.

The uniqueness of the solution of this CIP is an open and long standing
problem. In fact, uniqueness of a similar coefficient inverse problem can be
currently proven only in the case if the right hand side of equation (\ref%
{eq:helmholtz}) is a function which is not vanishing in $\overline{\Omega }.$
This can be done by the Bukhgeim-Klibanov method \cite{KlibanovBukhgeim81},
also see, e.g. \cite{BY,KlibanovTimonov04,Ksurvey} and references cited
therein for this method. Hence, for the computational purpose, we assume
below the uniqueness of our CISP.

In this last part of this section we want to briefly describe the travel
time $\tau (\mathbf{x})$ which is important in our analysis. The Riemannian
metric generated by the function $c(\mathbf{x})$ is: 
\begin{equation*}
d\tau (\mathbf{x})=\sqrt{c(\mathbf{x})}|d\mathbf{x}|,\quad |d\mathbf{x}|=%
\sqrt{(dx)^{2}+(dy)^{2}+(dz)^{2}}.
\end{equation*}%
For a fixed number $a>0$, consider the plane $P_{a}=\{(x,y,-a):x,y\in 
\mathbb{R}\}.$ We assume that $\Omega \subset \left\{ z>-a\right\} $ and
impose everywhere below the following condition on the function $c(\mathbf{x}%
)$:

\textbf{Regularity Assumption}. \emph{For any point }$\mathbf{x}\in \mathbb{R%
}^{3}$\emph{\ there exists a unique geodesic line }$\Gamma (\mathbf{x},a)$%
\emph{, with respect to the metric }$d\tau $\emph{, connecting }$\mathbf{x}$%
\emph{\ with the plane }$P_{a}$\emph{\ and perpendicular to }$P_{a}$\emph{\
near the intersection point.}

A sufficient condition of the regularity of geodesic lines is \cite{Rom}: 
\begin{equation*}
\sum\limits_{i,j=1}^{3}\frac{\partial ^{2}\ln (c(\mathbf{x}))}{\partial
x_{i}\partial x_{j}}\xi _{i}\xi _{j}\geq 0,\quad \text{for all }\mathbf{x}%
\in \overline{\Omega },\mathbf{\xi }\in \mathbb{R}^{3}.
\end{equation*}%
We introduce the travel time $\tau (\mathbf{x})$\ from the plane $P_{a}$ to
the point $\mathbf{x}$ as \cite{KlibanovRomanov16}\emph{\ }%
\begin{equation*}
\tau (\mathbf{x})=\int\limits_{\Gamma (\mathbf{x},a)}\sqrt{c\left( \mathbf{%
\xi }\right) }d\sigma .
\end{equation*}

\section{The Integro-Differential Equation}

\label{sec:3} %

In this section we reformulate our coefficient inverse problems as an
integro-differential equation, which is one of the main ingredients in our
reconstruction method. To this end, we first need a result on (high
frequency) asymptotic behavior of the total field $u(\mathbf{x},k)$ in \cite%
{KlibanovRomanov16}. It was shown in this cited paper that 
\begin{equation}
u(\mathbf{x},k)=A(\mathbf{x})e^{ik\left( \tau (\mathbf{x})-a\right) }\left[
1+s(\mathbf{x},k)\right] ,\quad \mathbf{x}\in \overline{\Omega }%
,k\rightarrow \infty ,  \label{eq:uasymptotics}
\end{equation}%
where $A(\mathbf{x})>0$ and $s(\mathbf{x},k)$ satisfies 
\begin{equation}
s(\mathbf{x,}k)=O\left( \frac{1}{k}\right) ,\quad \partial _{k}s(\mathbf{x,}%
k)=O\left( \frac{1}{k}\right) ,\quad \mathbf{x}\in \overline{\Omega }%
,k\rightarrow \infty .  \label{3.1}
\end{equation}%
Here $\tau (\mathbf{x})$ is the length of the geodesic line generated by the
function $c(\mathbf{x})$ in the Riemannian metric. Define 
\begin{equation}
w(\mathbf{x},k)=\frac{u(\mathbf{x},k)}{u_{i}(\mathbf{x},k)}.  \label{eq:w}
\end{equation}%
From (\ref{eq:uinc}), (\ref{eq:uasymptotics}) and (\ref{eq:w}), we have 
\begin{equation}
w(\mathbf{x},k)=A(\mathbf{x})e^{ik(\tau (\mathbf{x})-z-a)}[1+s(\mathbf{x}%
,k)],\quad \mathbf{x}\in \overline{\Omega },k\rightarrow \infty .
\label{eq:wasymptotics}
\end{equation}%
From (\ref{eq:uasymptotics}) and (\ref{eq:wasymptotics}), and for $\mathbf{x}%
\in \Omega $, $k\in \lbrack \underline{k},\overline{k}]$, we can uniquely
define the function $\log w(\mathbf{x},k)$ for sufficiently large values of $%
\underline{k}$ as 
\begin{equation}
\log w(\mathbf{x},k)=\ln A(\mathbf{x})+ik(\tau (\mathbf{x})-z-a)+%
\mathop{\displaystyle \sum }\limits_{n=1}^{\infty }\frac{\left( -1\right)
^{n-1}}{n}\left( s(\mathbf{x},k)\right) ^{n}.  \label{3.2}
\end{equation}%
It is clear that, with $\log w(\mathbf{x},k)$ defined as above, $\exp [\log
w(\mathbf{x},k)]$ equals to the right hand side of (\ref{eq:wasymptotics}).
Thus, we assume below that the number $\underline{k}$ is sufficiently large.

%

Now we are ready to derive the integro-differential equation. For $\mathbf{x}%
\in \Omega ,\,k\in \lbrack \underline{k},\overline{k}]$ we define the
function $v(\mathbf{x},k),$ 
\begin{equation}
v(\mathbf{x},k)=\frac{\log w(\mathbf{x},k)}{k^{2}}.  \label{eq:v}
\end{equation}%
Then 
\begin{equation}
\Delta v+k^{2}\nabla v\cdot \nabla v=-c(\mathbf{x}).  \label{3.20}
\end{equation}%
Setting $q(\mathbf{x},k)$ as 
\begin{equation}
q(\mathbf{x},k)=\partial _{k}v(\mathbf{x},k),  \label{eq:q}
\end{equation}%
we obtain 
\begin{equation}
v(\mathbf{x},k)=-\int\limits_{k}^{\overline{k}}q(\mathbf{x},\kappa )d\kappa
+V(\mathbf{x}).  \label{eq:vq}
\end{equation}%
Here we call $V(\mathbf{x})$ the tail function, 
\begin{equation}
V(\mathbf{x})=v(\mathbf{x},\overline{k}).  \label{3.30}
\end{equation}%
Combining (\ref{eq:helmholtz}), (\ref{eq:uinc}), (\ref{eq:coef}) and (\ref%
{eq:w}), we obtain 
\begin{equation}
\Delta w+k^{2}\beta w+2ik\frac{\partial w}{\partial z}=0.
\label{eq:intdiffw}
\end{equation}%
Taking into account (\ref{eq:v}), equation (\ref{eq:intdiffw}) becomes 
\begin{equation}
\Delta v+k^{2}\nabla v\cdot \nabla v+2ik\frac{\partial v}{\partial z}+\beta (%
\mathbf{x})=0.  \label{eq:intdiffv}
\end{equation}%
To eliminate the function $\beta (\mathbf{x})$ we differentiate (\ref%
{eq:intdiffv}) with respect to $k,$ 
\begin{equation}
\Delta q+2k\nabla v\cdot \left( k\nabla q+\nabla v\right) +2i\left( k\frac{%
\partial q}{\partial z}+\frac{\partial v}{\partial z}\right) =0.
\label{eq:intdiffq}
\end{equation}%
Substituting (\ref{eq:vq}) into (\ref{eq:intdiffq}) leads to the following
integro-differential equation

\begin{equation}
\begin{split}
L(q)=\Delta q+2k\left( \nabla V-\int_{k}^{\overline{k}}\nabla q(\kappa
)d\kappa \right) \cdot & \left( k\nabla (q+ V)-\int_{k}^{\overline{k}}\nabla
q\left( \kappa \right) d\kappa \right) \\
+2i& \left( kq_{z}+V_{z}-\int_{k}^{\overline{k}}q_{z}\left(\kappa \right)
d\kappa \right) =0.
\end{split}
\label{eq:intdiff}
\end{equation}%
This equation is complemented with the overdetermined boundary conditions: 
\begin{equation}
\begin{split}
q(\mathbf{x},k)=\phi _{0}(\mathbf{x},k),\quad q_{z}(\mathbf{x},k)& =\phi
_{1}(\mathbf{x},k),\quad \mathbf{x}\in \Gamma ,\,k\in \lbrack \underline{k},%
\overline{k}], \\
q(\mathbf{x},k)& =0,\quad \mathbf{x}\in \partial \Omega \setminus \Gamma
,\,k\in [\underline{k},\overline{k}],
\end{split}
\label{eq:intdiffbcs}
\end{equation}%
where the functions $\phi _{0}$ and $\phi _{1}$ are computed from the
functions $g_{0}$ and $g_{1}$ in (\ref{eq:cisp}), (\ref{eq:gz0}). The third
boundary condition (\ref{eq:intdiffbcs}) follows from (\ref{eq:uinc}), (\ref%
{eq:w}), (\ref{eq:v}) and (\ref{eq:q}).

Note that in (\ref{eq:intdiff}) we have two unknowns $q(\mathbf{x},k)$ and $%
V(\mathbf{x})$. Hence, we will solve the problem (\ref{eq:intdiff}), (\ref%
{eq:intdiffbcs}) using a predictor-corrector method. Here we find some
approximation of $V(\mathbf{x})$ first and use it as a predictor, and then
solve for $q(\mathbf{x},k)$. One can see that if certain approximations of $%
q(\mathbf{x},k)$\ and $V(\mathbf{x})$ are found, then an approximation for
the unknown coefficient $c( \mathbf{x}) $\ can be found via (\ref{eq:vq})
and (\ref{3.20}) for a certain value of $k\in [ \underline{k},\overline{k}] $%
. In our computations we use $k=\underline{k}$ \ for that value. Therefore,
we focus below on approximating functions $q(\mathbf{x},k)$, $V(\mathbf{x})$.

\section{Approximation of the tail function}

In this section we present a method for finding an approximation of the tail
function $V(\mathbf{x})$. We note that this method is different the one
studied in \cite{KlibanovKolesov17}.

It follows from (\ref{3.2}) and (\ref{3.30}) that there exists a function $p(%
\mathbf{x})$ such that 
\begin{equation}
v( \mathbf{x},k) =\frac{p( \mathbf{x}) }{k}+O\left( \frac{1}{k^{2}}\right)
,\quad q( \mathbf{x},k) =-\frac{p( \mathbf{x}) }{k^{2}}+O\left( \frac{1}{%
k^{3}}\right) ,\quad k\rightarrow \infty ,\,\mathbf{x}\in \Omega .
\label{eq:vqasympt}
\end{equation}%
For sufficiently large $\overline{k}$, we drop $O(1/\overline{k}^{2}) $ and $%
O( 1/\overline{k}^{3}) $ in (\ref{eq:vqasympt}) and set 
\begin{equation}
v( \mathbf{x},k) =\frac{p( \mathbf{x}) }{k},\quad q( \mathbf{x},k) =-\frac{%
p( \mathbf{x}) }{k^{2}},\quad k\geq \overline{k}, \mathbf{x}\in \Omega .
\label{eq:vqtail}
\end{equation}%
Next, substituting (\ref{eq:vqtail}) in (\ref{eq:intdiff}) and setting $k=%
\overline{k}$, we obtain 
\begin{equation}
\Delta V(\mathbf{x})=0,\quad \mathbf{x}\in \Omega .  \label{eq:tail}
\end{equation}%
This equation is supplemented by the following boundary conditions: 
\begin{equation}
V(\mathbf{x})=\psi _{0}(\mathbf{x}),\quad V_{z}(\mathbf{x})=\psi _{1}(%
\mathbf{x}),\quad \mathbf{x}\in \Gamma ,\quad V(\mathbf{x})=0,\quad \mathbf{x%
}\in \partial \Omega \setminus \Gamma ,  \label{eq:tailbcs}
\end{equation}%
where functions $\psi _{0}$ and $\psi _{1}$ are computed using (\ref{eq:cisp}%
) and (\ref{eq:gz0}). Boundary conditions (\ref{eq:tailbcs}) are
over-determined ones. Due to the approximate nature of (\ref{eq:vqtail}), we
have observed that the obvious approach of finding the function $V(\mathbf{x}%
)$ by dropping the second boundary condition (\ref{eq:tailbcs}) and solving
the resulting Dirichlet boundary value problem for Laplace equation (\ref%
{eq:tail}) with the boundary data (\ref{eq:tailbcs}) does not provide
satisfactory results. The same observation was made in \cite%
{KlibanovKolesov17} for the 1D case.

We show in section \ref{sec:4} how do we approximately solve the problem (%
\ref{eq:tail})--(\ref{eq:tailbcs}).

\section{Partial Finite Differences}

\label{sec:4}

\subsection{Grid points}

\label{sec:4.1}

We now write differential operators in (\ref{eq:intdiff}) and (\ref{eq:tail}%
) in finite differences with respect to $x,y$. Let the domain $\Omega
_{1}\subset \mathbb{R}^{2}$ be the orthogonal projection of the domain $%
\Omega \subset \mathbb{R}^{3}$ in (\ref{eq:2.1}) on the plane $\left\{
z=0\right\} ,$%
\begin{equation*}
\Omega _{1}=\left\{ \left( x,y,z\right) :|x|<b,|y|<b,z=0\right\} .
\end{equation*}%
Consider a finite difference grid in $\Omega _{1}$ with the uniform grid
step size $h.$ This grid consists of points $\left\{ \left(
x_{j},y_{s}\right) \right\} _{j,s=1}^{N_{h}}\subset \overline{\Omega }_{1}.$
Denote 
\begin{equation}
\Omega _{h}=\{( x_{j},y_{s},z): \left( x_{j},y_{s}\right)
_{j,s=1}^{N_{h}}\subset \overline{\Omega }_{1},z\in \left( -\xi ,d\right) \}
.  \label{400.1}
\end{equation}%
For every interior point $\left( x_{j},y_{s},z\right) \in \overline{\Omega }%
\setminus \partial \Omega $ four neighboring points are:%
\begin{eqnarray*}
( x_{j+1},y_{s},z) =( x_{j}+h,y_{s},z), \quad ( x_{j-1},y_{s},z) =(
x_{j}-h,y_{s},z) , \\
( x_{j},y_{s+1},z) =( x_{j},y_{s}+h,z), \quad ( x_{j},y_{s-1},z) =(
x_{j},y_{s}-h,z) .
\end{eqnarray*}%
The corresponding Laplace operator written in partial finite differences is%
\begin{equation}
\Delta ^{h}u=u_{zz}+u_{xx}^{h}+u_{yy}^{h},  \label{400.2}
\end{equation}%
where $u_{xx}^{h}$ and $u_{yy}^{h}$ are finite difference analogs of
continuous derivatives $u_{xx}$ and $u_{yy},$%
\begin{equation}
u_{xx}^{h}( x_{j},y_{s},z) =\frac{u( x_{j}-h,y_{s},z) -2u( x_{j},y_{s},z)
+u( x_{j}+h,y_{s},z) }{h^{2}}  \label{400.3}
\end{equation}%
and similarly for $u_{yy}^{h}.$ Next, 
\begin{equation}
\nabla ^{h}u=( \partial _{x}^{h}u,\partial _{y}^{h}u,\partial _{z}u) ,
\label{400.4}
\end{equation}
where 
\begin{equation*}
\partial _{x}^{h}u( x_{j},y_{s},z) = \frac{u( x_{j}+h,y_{s},z) -u(
x_{j}-h,y_{s},z) }{2h} 
\end{equation*}
and similarly for $\partial _{y}^{h}u( x_{j},y_{s},z) .$

\subsection{Problems (\protect\ref{eq:intdiff})--(\protect\ref{eq:intdiffbcs}%
) and (\protect\ref{eq:tail})--(\protect\ref{eq:tailbcs}) in partial finite
differences}

\label{sec:4.2}

We now rewrite problem (\ref{eq:intdiff})--(\ref{eq:intdiffbcs}) in partial
finite differences. To this end, we keep in mind that only interior grid
points are involved in differential operators below. Using (\ref{400.1})--(%
\ref{400.4}), we obtain for $x\in \Omega _{h}$ 
\begin{equation}
\begin{split}
L^{h}(q)=\Delta ^{h}q+2k\left( \nabla ^{h}V-\int_{k}^{\overline{k}}\nabla
^{h}q( \kappa )d\kappa \right) \cdot & \left( k\nabla ^{h} (q+V)-\int_{k}^{%
\overline{k}}\nabla ^{h}q\left( \kappa \right) d\kappa \right) \\
+2i& \left( kq_{z}+V_{z}-\int_{k}^{\overline{k}}q_{z}\left(\kappa \right)
d\kappa \right) =0,
\end{split}
\label{400.5}
\end{equation}%
\begin{equation}
\begin{split}
q(\mathbf{x},k)=\phi _{0}(\mathbf{x},k),\quad q_{z}(\mathbf{x},k)& =\phi
_{1}(\mathbf{x},k),\quad \mathbf{x}\in \Gamma ,\,k\in \lbrack \underline{k},%
\overline{k}], \\
q(\mathbf{x},k)& =0,\quad \mathbf{x}\in \partial \Omega \setminus \Gamma
,\,k\in \lbrack \underline{k},\overline{k}].
\end{split}
\label{400.6}
\end{equation}

Similarly, problem (\ref{eq:tail})--(\ref{eq:tailbcs}) becomes 
\begin{equation}
\Delta ^{h}V(\mathbf{x})=0,\quad \mathbf{x}\in \Omega .  \label{400.7}
\end{equation}%
\begin{equation}
V(\mathbf{x})=\psi _{0}(\mathbf{x}),\quad V_{z}(\mathbf{x})=\psi _{1}(%
\mathbf{x}),\quad \mathbf{x}\in \Gamma ,\quad V(\mathbf{x})=0,\quad \mathbf{x%
}\in \partial \Omega \setminus \Gamma.  \label{400.8}
\end{equation}

\begin{remark}
\label{rem:1} \normalfont 
\mbox{}

\begin{enumerate}
\item From now on functions $q( \mathbf{x},k) $ and $V( \mathbf{x})$, and
other functions we consider are semidiscrete, i.e. they are defined on $%
\overline{\Omega }_{h}.$ This means that, e.g. $q( \mathbf{x},k) =\left\{ q(
x_{j},y_{s},z) \right\} _{j,s=1}^{N_{h}}$, $V( \mathbf{x}) =\left\{ V(
x_{j},y_{s},z) \right\} _{j,s=1}^{N_{h}},$ etc. Boundary conditions at $%
\partial \Omega $ for the functions $q$ and $V$ are also defined only on
grid points which belong to the boundary $\partial \Omega $.

\item Since the grid step size $h$ is not changing in our arrangement, we
will not indicate below for brevity the dependence of some parameters on $h$%
, although they do depend on $h$. Thus, for example below $C=C(\xi ,d)>0$
denotes different positive constants depending only on numbers $\xi $,$d$
and $h$.
\end{enumerate}
\end{remark}

\subsection{Some functional spaces}

\label{sec:4.3}

Denote by $\overline{z}$ the complex conjugate of $z\in \mathbb{C}$. It is
convenient for us to consider any complex valued function $U=\mathop{\rm Re}%
U+i\mathop{\rm
Im}U=U_{1}+iU_{2}$ as the 2D vector function $U=\left( U_{1},U_{2}\right) .$
Furthermore, each component $U_{j}$ of this vector function is, in turn,
another vector function defined on the above grid, $%
U_{j}=U_{j}(x_{j},y_{s},z,k).$ Hence, below any Banach space of complex
valued functions is actually the space of these real valued vector functions
with the well known definitions of norms and scalar products (if in Hilbert
spaces). For brevity we do not differentiate below between complex valued
functions and corresponding vector functions. These things are always clear
from the context.

We introduce the Hilbert spaces $H^{2,h}(\Omega _{h}),L_{2}^{h}\left( \Omega
_{h}\right) $ and $H_{n}^{h}$ of semidiscrete complex valued functions as%
\begin{align*}
& H^{2,h}(\Omega _{h})=\{f(x_{j},y_{s},z):\left\Vert f\right\Vert
_{H^{n,h}(\Omega
_{h})}^{2}=\sum\limits_{j,s=1}^{N_{h}}\sum\limits_{r=0}^{2}h^{2}\int%
\limits_{-\xi }^{d}\left\vert \partial _{z}^{r}f\left( x_{j},y_{s},z\right)
\right\vert ^{2}dz<\infty \}, \\
& L_{2}^{h}\left( \Omega _{h}\right) =\{f(x_{j},y_{s},z):\left\Vert
f\right\Vert _{L_{2}^{h}(\Omega
_{h})}^{2}=\sum\limits_{j,s=1}^{N_{h}}h^{2}\int\limits_{-\xi }^{d}\left\vert
f\left( x_{j},y_{s},z\right) \right\vert ^{2}dz<\infty \}, \\
& H_{n}^{h}=\{f(x_{j},y_{s},z,k):\left\Vert f\right\Vert
_{H_{n}^{h}}^{2}=\int\limits_{\underline{k}}^{\overline{k}}\left\Vert
f\left( \mathbf{x},k\right) \right\Vert _{H^{n,h}\left( \Omega _{h}\right)
}^{2}dk<\infty \},\quad n=2,3.
\end{align*}%
Denote $\left[ ,\right] $ the scalar product in the space $H^{2,h}(\Omega
_{h})$. We also define subspaces $H_{0}^{2,h}(\Omega _{h})\subset
H^{2,h}(\Omega _{h})$ and $H_{0,2}^{h}\subset H_{2}^{h}$ as%
\begin{align*}
& H_{0}^{2,h}(\Omega _{h})=\{f(x_{j},y_{s},z)\in H^{2,h}(\Omega _{h}):f(%
\mathbf{x})\mid _{\partial \Omega }=0,f_{z}(\mathbf{x})\mid _{\Gamma }=0\},
\\
& H_{0,2}^{h}=\{f(x_{j},y_{s},z,k)\in H_{2}^{h}:f(\mathbf{x},k)\mid
_{\partial \Omega }=0,f_{z}(\mathbf{x},k)\mid _{\Gamma }=0,\forall k\in
\lbrack \underline{k},\overline{k}]\}.
\end{align*}%
Note that since, for all $f\in H_{0,2}^{h}$, 
\begin{equation*}
f(x_{j},y_{s},z,k)=\int\limits_{-\xi }^{z}f_{z}(x_{j},y_{s},\rho ,k)d\rho
,\quad f_{z}(x_{j},y_{s},z,k)=\int\limits_{-\xi }^{z}f_{zz}(x_{j},y_{s},\rho
,k)d\rho ,
\end{equation*}%
then the norm in $H_{0,2}^{h}$ is equivalent with%
\begin{equation}
\left\Vert f(\mathbf{x})\right\Vert _{H_{0,2}^{h}\left( \Omega _{h}\right)
}^{2}=\sum\limits_{j,s=1}^{N_{h}}h^{2}\int\limits_{-\xi }^{d}\left\vert
\Delta ^{h}f(x_{j},y_{s},z)\right\vert ^{2}dz.  \label{100}
\end{equation}%
In addition, for $l=0,1$%
\begin{align*}
C^{l}(\overline{\Omega }_{h})& =\{f(x_{j},y_{s},z):\left\Vert f\right\Vert
_{C^{l}(\overline{\Omega }_{h})}=\max_{j,s}\left\Vert
f(x_{j},y_{s},z)\right\Vert _{C^{l}\left[ -\xi ,d\right] }<\infty \}, \\
C_{l}^{h}& =\{f(x_{j},y_{s},z,k):\left\Vert f\right\Vert
_{C_{l}^{h}}=\max_{k\in \lbrack \underline{k},\overline{k}]}\left\Vert
f(x_{j},y_{s},z,k)\right\Vert _{C^{l}(\overline{\Omega }_{h})}<\infty \}.
\end{align*}%
By embedding theorem $H^{2,h}(\Omega _{h})\subset C^{1}(\overline{\Omega }%
_{h}),H_{n}^{h}\subset C_{n-1}^{h}$ and 
\begin{align}
& \left\Vert f\right\Vert _{C^{1}(\overline{\Omega }_{h})}\leq C\left\Vert
f\right\Vert _{H^{2,h}(\Omega _{h})},\quad \text{for all }f\in
H^{2,h}(\Omega _{h}),\text{ }  \label{400.9} \\
& \left\Vert f\right\Vert _{C_{n-1}^{h}}\leq C\left\Vert f\right\Vert
_{H_{n}^{h}},\quad \text{for all }f\in H_{n}^{h}.  \label{400.10}
\end{align}

\section{Two Cost Functionals with CWFs}

\label{sec:6}

It is our computational experience for the 1D case \cite%
{BakushinskiiKlibanov17, KlibanovKolesov17, KlibanovThanh15} that one should
use for computations such a CWF which would be a simple one. A similar
conclusion can be found on page 1581 of \cite{Baud}. Thus, the CWF we use in
this paper is: 
\begin{equation}
\varphi _{\lambda }(z)=e^{-2\lambda z}.  \label{500.1}
\end{equation}

\subsection{Problem (\protect\ref{400.7})--(\protect\ref{400.8})}

\label{sec:6.1}

First, we present the cost functional for the solution of problem (\ref%
{400.7})--(\ref{400.8}) which is about the tail function. Non-zero boundary
conditions in (\ref{400.8}) are inconvenient for us. Hence, we assume that
there exists a function $Q(\mathbf{x})\in H^{2,h}(\Omega _{h})$ such that%
\begin{equation}
Q(\mathbf{x})=\psi _{0}(\mathbf{x}),\quad \partial _{z}Q(\mathbf{x})=\psi
_{1}(\mathbf{x}),\quad \mathbf{x}\in \Gamma ;\quad Q(\mathbf{x})=0,\quad 
\mathbf{x}\in \partial \Omega \setminus \Gamma ,  \label{500.2}
\end{equation}%
Define 
\begin{equation}
W(\mathbf{x})=V(\mathbf{x})-Q(\mathbf{x})\in H_{0}^{2,h}(\Omega _{h}).
\label{500.20}
\end{equation}%
Hence, we consider the following minimization problem:

\textbf{Minimization Problem 1}. \emph{For $W\in H_{0}^{2,h}(\Omega _{h})$,
minimize the functional }$I_{\mu }(W),$ 
\begin{equation}
I_{\mu }\left( W\right) =e^{2\mu
d}\sum\limits_{j,s=1}^{N_{h}}h^{2}\int\limits_{-\xi }^{d}\left\vert (\Delta
^{h}W+\Delta ^{h}Q)(x_{j},y_{s},z)\right\vert ^{2}\varphi _{\mu }(z)dz.
\label{500.3}
\end{equation}

The multiplier $e^{2\mu d}$ is introduced here to ensure that $e^{2\mu
d}\min_{\left[ -\xi ,d\right] }\varphi _{\mu }(z)=1.$

\begin{remark}
Since the operator $\Delta ^{h}$ is linear, then, in principle at least, one
can apply straightforwardly the quasi-reversibility method to find an
approximate solution of the problem $\Delta W+\Delta Q=0$ for $W\in
H_{0}^{2,h}(\Omega _{h})$ \cite{KQR}. This means that one can use $\lambda =0
$ in (\ref{500.3}). However, it was observed in \cite{BakushinskiiKlibanov17}
that the involvement of the CWF like in (\ref{500.3}) leads to a better
solution accuracy.
\end{remark}

We now follow the classical Tikhonov regularization concept \cite{Bak,T}. By
this concept, we should assume that there exists an exact solution $V_{\ast
}(\mathbf{x})$ of the problem (\ref{400.7})--(\ref{400.8}) with the
noiseless data $\psi _{0\ast }(\mathbf{x}),\psi _{1\ast }(\mathbf{x}).$
Below the subscript \textquotedblleft $\ast $" is related only to the exact
solution. In fact, however, the data $\psi _{0}(\mathbf{x})$ and $\psi _{1}(%
\mathbf{x})$ contain noise. Let $\delta \in \left( 0,1\right) $ be the level
of noise in the data $\psi _{0}(\mathbf{x})$ and $\psi _{1}(\mathbf{x})$.
Again, following the same concept, we should assume that the number $\delta
\in \left( 0,1\right) $ is sufficiently small. Assume that there exists the
function $Q_{\ast }(\mathbf{x})\in H^{2,h}(\Omega _{h})$ such that

\begin{equation}
Q_{\ast }(\mathbf{x})=\psi _{0\ast }(\mathbf{x}),\quad \partial _{z}Q_{\ast
}(\mathbf{x})=\psi _{1\ast }(\mathbf{x}),\quad \mathbf{x}\in \Gamma ;\quad
Q_{\ast }(\mathbf{x})=0,\quad \mathbf{x}\in \partial \Omega \setminus \Gamma
,  \label{3.5}
\end{equation}%
\begin{equation}
\left\Vert Q-Q_{\ast }\right\Vert _{H^{2,h}\left( \Omega _{h}\right)
}<\delta ,  \label{3.6}
\end{equation}%
where $Q$ is defined in (\ref{500.2}). We will choose in Theorem 7.2 of
section \ref{sec:6} a certain dependence $\mu =\mu \left( \delta \right) $
of the parameters $\mu $ on the noise level $\delta $. Denote $W_{\mu
(\delta )}(\mathbf{x})=W_{\min }(\mathbf{x})$ the unique minimizer of the
functional $I_{\mu (\delta )}(W)$ (Theorem 7.2) and by (\ref{500.20}) let 
\begin{equation}
V_{\mu (\delta )}(\mathbf{x})=W_{\mu (\delta )}(\mathbf{x})+Q(\mathbf{x}%
)=W_{\min }(\mathbf{x})+Q(\mathbf{x}).  \label{3.7}
\end{equation}

\subsection{Problem (\protect\ref{eq:intdiff})--(\protect\ref{eq:intdiffbcs})%
}

\label{sec:6.2}

Suppose that there exists a function $F(\mathbf{x},k)\in H_{3}^{h}$ such
that (see (\ref{eq:intdiffbcs})): 
\begin{equation}
F(\mathbf{x},k)=\phi _{0}(\mathbf{x},k),\quad F_{z}(\mathbf{x},k)=\phi _{1}(%
\mathbf{x},k),\quad \mathbf{x}\in {\Gamma },\quad F(\mathbf{x},k)=0,\quad 
\mathbf{x}\in {\partial \Omega \setminus \Gamma }.  \label{3.8}
\end{equation}%
Also, assume that there exists an exact solution $c_{\ast }(\mathbf{x})$ of
our CIP satisfying the above conditions imposed on the coefficient $c(%
\mathbf{x})$ and generating the noiseless boundary data $\phi _{0,\ast }$
and $\phi _{1,\ast }$ in (\ref{eq:intdiffbcs}). Also, assume that there
exists the function $F_{\ast }(\mathbf{x},k)\in H_{3}^{h}$ satisfying the
following analog of boundary conditions (\ref{3.8}):%
\begin{equation}
F_{\ast }(\mathbf{x},k)=\phi _{0,\ast }(\mathbf{x},k),\,\partial _{z}F_{\ast
}(\mathbf{x},k)=\phi _{1,\ast }(\mathbf{x},k),\,\mathbf{x}\in {\Gamma }%
,\,F_{\ast }(\mathbf{x},k)=0,\,\mathbf{x}\in {\partial \Omega \setminus
\Gamma .}  \label{3.88}
\end{equation}%
We assume that%
\begin{equation}
\left\Vert F-F_{\ast }\right\Vert _{H_{3}^{h}}<\delta .  \label{3.80}
\end{equation}%
Let $q_{\ast }\in H_{2}^{h}$ be the function $q$ generated by the exact
coefficient $c_{\ast }(\mathbf{x})$. We define functions $p$ and $p_{\ast }$
as 
\begin{equation}
p(\mathbf{x},k)=q(\mathbf{x},k)-F(\mathbf{x},k),\quad p_{\ast }(\mathbf{x}%
,k)=q_{\ast }(\mathbf{x},k)-F_{\ast }(\mathbf{x},k).  \label{3.9}
\end{equation}%
Hence, the functions $p$,$p_{\ast }\in $ $H_{0,2}^{h}.$ Let $R>0$ be an
arbitrary number. Consider the ball $B(R)\subset H_{0,2}^{h}$ of the radius $%
R$, 
\begin{equation}
B(R)=\{r\in H_{0,2}^{h}:\left\Vert r\right\Vert _{H_{2}^{h}}<R\}.
\label{3.10}
\end{equation}

Using the integro-differential equation (\ref{eq:intdiff}), boundary
conditions (\ref{eq:intdiffbcs}) for it, (\ref{3.8}), (\ref{3.88}) and (\ref%
{3.9}), we construct our cost functional $J_{\lambda }(p)$ with the CWF (\ref%
{500.1}) in it as: 
\begin{equation}
J_{\lambda }(p)=e^{2\lambda d}\sum\limits_{j,s=1}^{N_{h}}h^{2}\int\limits_{%
\underline{k}}^{\overline{k}}\int\limits_{-\xi
}^{d}|L^{h}(p+F)(x_{j},y_{s},z,\kappa )|^{2}\varphi _{\lambda }(z)dzd\kappa
,\quad p\in \overline{B(R)},  \label{500.4}
\end{equation}%
where the tail function in $L^{h}$ is defined in (\ref{3.7}). Similarly with
(\ref{500.3}), the multiplier $e^{2\lambda d}$ is introduced to balance two
terms in the right hand side of (\ref{500.4}). We consider the following
minimization problem:

\noindent \textbf{Minimization Problem 2}. \emph{Minimize the functional }$%
J_{\lambda }(p)$\emph{\ on the set }$p\in \overline{B\left( R\right) }.$

\section{Carleman Estimate and Global Strict Convexity}

\label{sec:7}

In this section we formulate theorems about the minimization problems 1 and
2 of section~\ref{sec:6}. First, we are concerned with the Carleman estimate
with the CWF (\ref{500.1}).

\textbf{Theorem 7.1 }(Carleman estimate).\emph{\ For }$\lambda >0$\emph{\
let }%
\begin{equation*}
B_{h}(u,\lambda )=\sum\limits_{j,s=1}^{M_{h}}h^{2}\int\limits_{-\xi
}^{d}|\Delta ^{h}u(x_{j},y_{s},z)|^{2}\varphi _{\lambda }(z)dz.
\end{equation*}%
\emph{Then there exists a sufficiently large number }$\lambda _{0}=\lambda
_{0}(\xi ,d)>1$\emph{\ such that for all $\lambda \geq \lambda _{0}$ the
following estimate is valid for all functions }$u\in H_{0}^{2,h}(\Omega _{h})
$\emph{\ } 
\begin{align}
B_{h}(u,\lambda )& \geq C\sum\limits_{j,s=1}^{M_{h}}h^{2}\int\limits_{-\xi
}^{d}\left\vert u_{zz}\left( x_{j},y_{s},z\right) \right\vert ^{2}\varphi
_{\lambda }\left( z\right) dz+C\lambda
\sum\limits_{j,s=1}^{M_{h}}h^{2}\int\limits_{-\xi }^{d}\left\vert
u_{z}\left( x_{j},y_{s},z\right) \right\vert ^{2}\varphi _{\lambda }(z)dz 
\notag \\
& +C\lambda ^{3}\sum\limits_{j,s=1}^{M_{h}}h^{2}\int\limits_{-\xi
}^{d}\left\vert u\left( x_{j},y_{s},z\right) \right\vert ^{2}\varphi
_{\lambda }(z)dz.  \label{6.0}
\end{align}

\textit{Proof.} Recall that we do not indicate the dependence of neither
constants $C$ nor other constants on $h$ (second item in Remarks 5.1). Since 
$|v|^{2}=(\mathop{\rm Re}v)^{2}+(\mathop{\rm Im}v)^{2}$ for $v\in \mathbb{C},
$ then it is sufficient to prove estimate (\ref{6.0}) for real valued
functions $u\in H_{0}^{2,h}(\Omega _{h})$\emph{. }The following Carleman
estimate was proven in lemma 3.1 of \cite{KlibanovKolesov17} for all real
valued functions $w(z)\in H^{2}(-\xi ,d)$ such that $w(-\xi )=w^{\prime
}(-\xi )=0$: 
\begin{equation}
\int\limits_{-\xi }^{d}(w^{\prime \prime })^{2}\varphi _{\lambda }(z)dz\geq 
\widetilde{C}\int\limits_{-\xi }^{d}(w^{\prime \prime })^{2}\varphi
_{\lambda }(z)dz+\widetilde{C}\lambda \int\limits_{-\xi }^{d}(w^{\prime
})^{2}\varphi _{\lambda }(z)dz+\widetilde{C}\lambda ^{3}\int\limits_{-\xi
}^{d}w^{2}\varphi _{\lambda }(z)dz,  \label{7.5}
\end{equation}%
for all $\lambda \geq \lambda _{0}(\xi ,d),$ where the number $\widetilde{C}=%
\widetilde{C}(\xi ,d)>0$ depends only on $\xi $ and $d.$ Next, it follows
from (\ref{400.2}) and (\ref{400.3}) that 
\begin{equation*}
B_{h}(u,\lambda )=\sum\limits_{j,s=1}^{M_{h}}h^{2}\int\limits_{-\xi
}^{d}[(u_{zz}+u_{xx}^{h}+u_{yy}^{h})(x_{j},y_{s},z)]^{2}\varphi _{\lambda
}(z)dz
\end{equation*}%
\begin{equation*}
\geq \frac{1}{2}\sum\limits_{j,s=1}^{M_{h}}h^{2}\int\limits_{-\xi
}^{d}[u_{zz}(x_{j},y_{s},z)]^{2}\varphi _{\lambda
}(z)dz-C\sum\limits_{j,s=1}^{M_{h}}\int\limits_{-\xi
}^{d}[u(x_{j},y_{s},z)]^{2}\varphi _{\lambda }(z)dz.
\end{equation*}%
Hence, using (\ref{7.5}), we obtain%
\begin{align}
B_{h}(u,\lambda )& \geq C\sum\limits_{j,s=1}^{M_{h}}h^{2}\int\limits_{-\xi
}^{d}[u_{zz}(x_{j},y_{s},z)]^{2}\varphi _{\lambda }(z)dz+C\lambda
\sum\limits_{j,s=1}^{M_{h}}h^{2}\int\limits_{-\xi
}^{d}[u_{z}(x_{j},y_{s},z)]^{2}\varphi _{\lambda }(z)dz  \notag \\
& +C\lambda ^{3}\sum\limits_{j,s=1}^{M_{h}}h^{2}\int\limits_{-\xi
}^{d}[u(x_{j},y_{s},z)]^{2}\varphi _{\lambda
}(z)dz-C\sum\limits_{j,s=1}^{M_{h}}\int\limits_{-\xi
}^{d}[u(x_{j},y_{s},z)]^{2}\varphi _{\lambda }(z)dz.  \label{7.6}
\end{align}%
Now choosing $\lambda _{0}$ so large that $\lambda _{0}^{3}h^{2}>C/2$, we
obtain from (\ref{7.6}) the target estimate (\ref{6.0}) for all $\lambda
>\lambda _{0}$. $\square $

The next theorem is about the functional $I_{\mu }(W)$ in (\ref{500.3}).

\textbf{Theorem 7.2.} \emph{Assume that there exists a function }$Q\in
H^{2,h}(\Omega _{h})$\emph{\ satisfying conditions (\ref{3.5}). Introduce
the function }$W\in H_{0}^{2,h}(\Omega _{h})$\emph{\ via (\ref{500.20}).
Then for each }$\mu >0$\emph{\ there exists unique minimizer }$W_{\mu }\in
H_{0}^{2,h}(\Omega _{h})$\emph{\ of the functional (\ref{500.3}). Suppose
now that there exists an exact solution }$V_{\ast }\in H^{2,h}(\Omega _{h})$%
\emph{\ of equation (\ref{eq:tail}) with the boundary data }$\psi _{0\ast
}(x)$ and $\psi _{1\ast }(x)$\emph{\ in (\ref{eq:tailbcs}). Also, assume
that there exists a function }$Q_{\ast }\in H^{2,h}(\Omega _{h})$\emph{\
satisfying conditions (\ref{3.5}) and such that inequality (\ref{3.6})
holds, where }$\delta \in (0,1)$\emph{\ is the noise level in the data. Let }%
$\lambda _{0}>0$\emph{\ be the number of Theorem 6.1. Choose a number }$%
\delta _{0}\in (0,e^{-2\left( d+\xi \right) \lambda _{0}})$. \emph{For any} $%
\delta \in \left( 0,\delta _{0}\right) $\emph{\ let }%
\begin{equation}
\mu =\mu (\delta )=\ln (\delta ^{-1/(2\left( d+\xi \right) )}).
\label{6.100}
\end{equation}%
\emph{\ Let the function }$V_{\mu (\delta )}(\mathbf{x})$ \emph{be defined
via (\ref{3.7})}. \emph{Then the following convergence estimate of} $V_{\mu
(\delta )}(\mathbf{x})$\emph{\ to the exact solution }$V_{\ast }(\mathbf{x})$
\emph{holds as }$\delta \rightarrow 0$ 
\begin{equation}
\left\Vert V_{\mu \left( \delta \right) }-V_{\ast }\right\Vert
_{H^{2,h}(\Omega _{h})}\leq C\sqrt{\delta }.  \label{6.2}
\end{equation}%
\emph{In addition,} $V_{\mu (\delta )}\in C^{1}(\overline{\Omega }_{h})$ 
\emph{\ and\ } 
\begin{equation}
C\left\Vert \nabla V_{\mu (\delta )}\right\Vert _{C^{1}(\overline{\Omega }%
_{h})}\leq \left\Vert V_{\mu (\delta )}\right\Vert _{H^{2,h}(\Omega
_{h})}\leq C\left[ 1+\left\Vert V_{\ast }\right\Vert _{H^{2,h}(\Omega _{h})}%
\right] .  \label{6.3}
\end{equation}

\textit{Proof.} It follows from (\ref{500.3}) and the variational principle
that the vector function $W_{\min }=(W_{1,\min },W_{2,\min })\in
H_{0}^{2,h}(\Omega _{h})$ is a minimizer of the functional $I_{\mu ,\alpha
}(W)$ if and only if 
\begin{align}
& e^{2\mu d}\sum\limits_{j,s=1}^{M_{h}}h^{2}\int\limits_{-\xi }^{d}(\Delta
^{h}W_{1,\min }\Delta ^{h}r_{1}+\Delta ^{h}W_{2,\min }\Delta
^{h}r_{2})(x_{j},y_{s},z)\varphi _{\mu }(z)dz  \notag \\
& =-e^{2\mu d}\sum\limits_{j,s=1}^{M_{h}}h^{2}\int\limits_{-\xi }^{d}(\Delta
^{h}Q_{1}\Delta ^{h}r_{1}+\Delta ^{h}Q_{2}\Delta
^{h}r_{2})(x_{j},y_{s},z)\varphi _{\mu }(z)dz,  \label{7.2}
\end{align}%
for all $r=(r_{1},r_{2})\in H_{0}^{2,h}(\Omega _{h})$. For any vector
function $P=(P_{1},P_{2})\in H_{0}^{2,h}(\Omega _{h})$ consider the
expression in the left hand side of (\ref{7.2}) in which the vector function 
$(W_{1,\min },W_{2,\min })$ is replaced with $(P_{1},P_{2})$. Then (\ref{100}%
) implies that this expression defines a new scalar product $\{P,r\}$ in the
space $H_{0}^{2,h}(\Omega _{h}),$\ and the corresponding norm $\{P,P\}^{1/2}$
is equivalent to the norm in the space $H^{2,h}(\Omega _{h})$. Next, for all 
$r=(r_{1},r_{2})\in H_{0}^{2,h}(\Omega _{h})$, we have 
\begin{align*}
\left\vert -e^{2\mu d}\sum\limits_{j,s=1}^{M_{h}}h^{2}\int\limits_{-\xi
}^{d}(\Delta ^{h}Q_{1}\Delta ^{h}r_{1}+\Delta ^{h}Q_{2}\Delta
^{h}r_{2})(x_{j},y_{s},z)\varphi _{\mu }(z)dz\right\vert & \leq D\left\Vert
Q\right\Vert _{H^{2,h}(\Omega _{h})}\left\Vert r\right\Vert _{H^{2,h}(\Omega
_{h})} \\
& \leq D_{1}\sqrt{\left\{ Q,Q\right\} }\sqrt{\left\{ r,r\right\} },
\end{align*}%
where the constants $D,D_{1}$ do not depend on $Q$ and $r$. Hence, by Riesz
theorem there exists unique vector function $\widehat{Q}=\left( \widehat{Q}%
_{1},\widehat{Q}_{2}\right) =\widehat{Q}(Q)\in H_{0}^{2,h}(\Omega _{h})$
such that 
\begin{equation*}
\{\widehat{Q},r\}=-e^{2\mu
d}\sum\limits_{j,s=1}^{M_{h}}h^{2}\int\limits_{-\xi }^{d}(\Delta
^{h}Q_{1}\Delta ^{h}r_{1}+\Delta ^{h}Q_{2}\Delta
^{h}r_{2})(x_{j},y_{s},z)\varphi _{\mu }(z)dz,
\end{equation*}%
for all $r=\left( r_{1},r_{2}\right) \in H_{0}^{2,h}(\Omega _{h})$. Hence,
by (\ref{7.2}) $\{W_{\min },r\}=\{\widehat{Q},r\}$, $\forall r\in
H_{0}^{2,h}(\Omega _{h})$. This implies that $W_{\min }=\widehat{Q}$. Thus,
both existence and uniqueness of the minimizer of the functional $I_{\mu }(W)
$ are established.

We now prove convergence estimate (\ref{6.2}). Let $W_{\ast }=\left( W_{\ast
,1},W_{\ast ,2}\right) =V_{\ast }-Q_{\ast }.$ Then $W_{\ast }\in
H_{0}^{2,h}(\Omega _{h}).$ Denote $\widetilde{W}=W_{\min }-W_{\ast }$ and $%
\widetilde{Q}=Q-Q_{\ast }$. Since 
\begin{align}
& e^{2\mu d}\sum\limits_{j,s=1}^{M_{h}}h^{2}\int\limits_{-\xi }^{d}(\Delta
^{h}W_{\ast ,1}\Delta ^{h}r_{1}+\Delta ^{h}W_{\ast ,2}\Delta
^{h}r_{2})(x_{j},y_{s},z)\varphi _{\mu }(z)dz  \notag \\
& =-e^{2\mu d}\sum\limits_{j,s=1}^{M_{h}}h^{2}\int\limits_{-\xi }^{d}(\Delta
^{h}Q_{1}^{\ast }\Delta ^{h}r_{1}+\Delta ^{h}Q_{2}^{\ast }\Delta
^{h}r_{2})(x_{j},y_{s},z)\varphi _{\mu }(z)dz,  \label{7.3}
\end{align}%
then subtracting (\ref{7.3}) from (\ref{7.2}) and setting $r=\widetilde{W},$
we obtain%
\begin{align*}
& e^{2\mu d}\sum\limits_{j,s=1}^{M_{h}}h^{2}\int\limits_{-\xi }^{d}(\Delta
^{h}\widetilde{W})^{2}\varphi _{\mu }(z)dz \\
& =-e^{2\mu d}\sum\limits_{j,s=1}^{M_{h}}h^{2}\int\limits_{-\xi }^{d}(\Delta
^{h}\widetilde{Q}_{1}\Delta ^{h}\widetilde{W}_{1}+\Delta ^{h}\widetilde{Q}%
_{2}\Delta ^{h}\widetilde{W}_{2})\varphi _{\mu }(z)dz.
\end{align*}%
Using the Cauchy-Schwarz inequality, taking into account (\ref{3.6}) and (%
\ref{6.100}), we obtain%
\begin{equation}
e^{2\mu d}\sum\limits_{j,s=1}^{M_{h}}h^{2}\int\limits_{-\xi }^{d}(\Delta ^{h}%
\widetilde{W}(x_{j},y_{s},z))^{2}\varphi _{\mu }(z)dz\leq Ce^{2\mu \left(
d+\xi \right) }\delta ^{2}.  \label{5.4}
\end{equation}%
By (\ref{6.100}) $\mu =\ln (\delta ^{-1/(2\left( d+\xi \right) )})$, which
implies $e^{2\mu \left( d+\xi \right) }\delta ^{2}=\delta .$ Hence, (\ref%
{5.4}) implies that%
\begin{equation}
e^{2\mu d}\sum\limits_{j,s=1}^{M_{h}}h^{2}\int\limits_{-\xi }^{d}(\Delta ^{h}%
\widetilde{W}(x_{j},y_{s},z))^{2}\varphi _{\mu }\left( z\right) dz\text{ }%
\mathbf{\leq }C\delta .  \label{5.50}
\end{equation}%
We now apply Theorem 7.1 to the left hand side of (\ref{5.50}). We obtain
for all $\mu \geq \mu _{0}$%
\begin{align}
& \delta \geq Ce^{2\mu d}\sum\limits_{j,s=1}^{M_{h}}h^{2}\int\limits_{-\xi
}^{d}(\widetilde{W}_{zz}(x_{j},y_{s},z))^{2}\varphi _{\mu }(z)dz  \notag \\
& +Ce^{2\mu d}\left[ \mu \sum\limits_{j,s=1}^{M_{h}}h^{2}\int\limits_{-\xi
}^{d}(\widetilde{W}_{z}(x_{j},y_{s},z))^{2}\varphi _{\mu }(z)dz+\mu
^{3}\sum\limits_{j,s=1}^{M_{h}}h^{2}\int\limits_{-\xi }^{d}(\widetilde{W}%
(x_{j},y_{s},z))^{2}\varphi _{\mu }(z)dz\right] .  \label{5.51}
\end{align}%
Since $e^{2\mu d}\varphi _{\mu }(z)\geq e^{2\mu d}\varphi _{\mu }(d)=1,$
then (\ref{5.51}) implies that 
\begin{equation*}
\Vert \widetilde{W}\Vert _{H^{2,h}(\Omega _{h})}\leq C\sqrt{\delta }.
\end{equation*}%
Hence, (\ref{3.6}), (\ref{3.7}) and triangle inequality imply that%
\begin{equation}
\left\Vert V_{\mu (\delta )}-V_{\ast }\right\Vert _{H^{2,h}(\Omega
_{h})}=\left\Vert \widetilde{W}+\left( Q-Q_{\ast }\right) \right\Vert
_{H^{2,h}(\Omega _{h})}\leq C\sqrt{\delta },  \label{7.4}
\end{equation}%
which proves (\ref{6.2}). Next, by (\ref{6.2}) and triangle inequality imply
the right estimate (\ref{6.3}). The left estimate (\ref{6.3}) follows from (%
\ref{400.9}). $\square $   

The main analytical result of this paper is Theorem 7.3.

\textbf{Theorem 7.3} (globally strict convexity). \emph{Assume that
conditions of Theorem 6.2 hold.} \emph{Let }$\lambda _{1}\geq \lambda _{0}$%
\emph{\ be the number defined below in the formulation of this theorem.
Assume that there exist functions }$F(\mathbf{x},k)$,$F_{\ast }(\mathbf{x}%
,k)\in H_{3}^{h}$ \emph{satisfying conditions (\ref{3.8})--(\ref{3.80}),
where }$\delta \in \left( 0,\delta _{1}\right) $ \emph{and} $\delta _{1}\in
(0,e^{-2\left( d+\xi \right) \lambda _{1}})$. \emph{Set in (\ref{500.4}) }$%
V=V_{\mu (\delta )},$\emph{\ where the function }$V_{\mu \left( \delta
\right) }$\emph{\ is defined in Theorem 7.2. First, the functional }$%
J_{\lambda }\left( p\right) $\emph{\ has the Frech\'{e}t derivative }$%
J_{\lambda }^{\prime }\left( p\right) \in H_{0,2}^{h}$\emph{\ at any point }$%
p\in H_{0,2}^{h}.$\emph{\ Second, there exist numbers } 
\begin{align*}
\lambda _{1}& =\lambda _{1}(\Omega _{h},R,\left\Vert F_{\ast }\right\Vert
_{H_{3}^{h}},\left\Vert V_{\ast }\right\Vert _{H^{2,h}(\Omega _{h})},%
\underline{k},\overline{k})\geq \lambda _{0}, \\
C_{1}& =C_{1}(\Omega _{h},R,\left\Vert F_{\ast }\right\Vert
_{H_{3}^{h}},\left\Vert V_{\ast }\right\Vert _{H^{2,h}(\Omega _{h})},%
\underline{k},\overline{k})>0,
\end{align*}%
\emph{depending only on listed parameters,} \emph{such that for any }$%
\lambda \geq \lambda _{1}$\emph{\ the functional }$J_{\lambda }(p)$\emph{\
is strictly convex on }$\overline{B(R)}.$\emph{\ In other words, the
following estimate holds:\ }%
\begin{equation}
J_{\lambda }(p+r)-J_{\lambda }(p)-J_{\lambda }^{\prime }(p)(r)\geq
C_{1}\left\Vert r\right\Vert _{H_{2}^{h}}^{2},\quad \text{for all }p,p+r\in 
\overline{B(R)}.  \label{6.4}
\end{equation}

\textit{Proof.} In this proof $C_{1}=C_{1}(\Omega _{h},R,\left\Vert F_{\ast
}\right\Vert _{H_{3}^{h}},\left\Vert V_{\ast }\right\Vert _{H^{2,h}(\Omega
_{h})},\underline{k},\overline{k})>0$ denotes different positive constants.
In addition, in this proof we denote for brevity $V(\mathbf{x})=V_{\mu
(\delta )}(\mathbf{x})$ and also sometimes we do not indicate the
dependencies on $(x_{j},y_{s},z,k)$. Note that (\ref{400.9}), (\ref{400.10}%
), (\ref{3.80})-(\ref{3.10}) and (\ref{6.3}) imply that%
\begin{equation}
\left\Vert \nabla V\right\Vert _{C^{1}(\overline{\Omega }_{h})}\leq
C_{1},\quad \left\Vert F\right\Vert _{C_{2}^{h}}\leq C_{1},  \label{5.52}
\end{equation}%
\begin{equation}
\Vert \nabla ^{h}r\Vert _{C_{1}^{h}}\leq C_{1}.  \label{5.53}
\end{equation}

Consider an arbitrary vector function $p=\left( p_{1},p_{2}\right) \in 
\overline{B(R)}$ and an arbitrary function $r=(r_{1},r_{2})\in H_{0,2}^{h}$
such that $p+r\in \overline{B(R)}$. By (\ref{500.4}) we need to consider $A$%
, where 
\begin{equation}
A=|L^{h}(p+r+F)|^{2}-|L^{h}(p+F)|^{2}.  \label{5.54}
\end{equation}%
First, we will single out such a part of $A$, which is linear with respect
to $r$. This will lead us to the Frech\'{e}t derivative $J_{\lambda
}^{\prime }.$ Next, we will single out $|\Delta ^{h}r|^{2}.$ Based on this,
we will apply the Carleman estimate of Theorem \ref{thm:1}. For all $%
z_{1},z_{2}\in \mathbb{C}$, we have

\begin{equation}
\left\vert z_{1}\right\vert ^{2}-\left\vert z_{2}\right\vert ^{2}=\left(
z_{1}-z_{2}\right) \overline{z}_{1}+\left( \overline{z}_{1}-\overline{z}%
_{2}\right) z_{2}.  \label{5.6}
\end{equation}%
Denote%
\begin{equation}
z_{1}=L^{h}(p+r+F),\text{ }z_{2}=L^{h}(p+F).  \label{5.7}
\end{equation}%
Then 
\begin{equation}
A_{1}=\left( z_{1}-z_{2}\right) \overline{z}_{1},\text{ }A_{2}=\left( 
\overline{z}_{1}-\overline{z}_{2}\right) z_{2},\quad A=A_{1}+A_{2}.
\label{5.81}
\end{equation}%
Using (\ref{400.5}), (\ref{3.9}) and (\ref{5.7}), we obtain%
\begin{equation}
\begin{gathered} z_{1}-z_{2} =\Delta ^{h}r-2k^{2}\nabla ^{h}r\cdot
\Bigg[\nabla ^{h}V-\int\limits_{k}^{\overline{k}}\nabla ^{h} (p+ F)d\kappa
\Bigg] +2i\Bigg[r_{z}-\int\limits_{k}^{\overline{k}}r_{z}d\kappa \Bigg]. \\
+2k\int\limits_{k}^{\overline{k}}\nabla ^{h}r \, d\kappa \cdot \Bigg[2\nabla
^{h}V-2\int\limits_{k}^{\overline{k}}\nabla ^{h} (p+ F)d\kappa +k \nabla
^{h} (p+ F)\Bigg] \end{gathered}  \label{5.9}
\end{equation}%
Next, 
\begin{equation*}
\begin{gathered} \overline{z}_{1}=\Delta ^{h} ( \overline{r}+ \overline{p}+
\overline{F}) -2k\Bigg [\nabla
^{h}\overline{V}-\int\limits_{k}^{\overline{k}}\nabla ^{h} ( \overline{r}+
\overline{p_{1}}+ \overline{F})d\kappa \Bigg]\cdot \Bigg[k\nabla ^{h}
(\overline{r}+\overline{p_{1}}+ \overline{F}) \\ +\nabla
^{h}\overline{V}-\int\limits_{k}^{\overline{k}}\nabla ^{h}(\overline{r}+
\overline{p}+ \overline{F})d\kappa \Bigg] -2i\Bigg[k\left(
\overline{r_{z}}+\overline{p_{z}}+\overline{F_{z}}\right)
+\overline{V_{z}}-\int\limits_{k}^{\overline{k}}\left(
\overline{r_{z}}+\overline{p_{z}}+\overline{F_{z}}\right) d\kappa \Bigg].
\end{gathered}
\end{equation*}%
Hence,%
\begin{equation}
A_{1}=(z_{1}-z_{2})\overline{z}_{1}=|\Delta
^{h}r|^{2}+Z_{l,1}(r,k)+Z_{1}(r,k),  \label{5.10}
\end{equation}%
where $Z_{l,1}(h,k)$ is linear with respect to the vector function $%
r=(r_{1},r_{2}),$%
\begin{equation}
\begin{gathered} Z_{l,1}\left( r,k\right) =\Delta ^{h}r\cdot Y_{1}+\nabla
^{h}r\nabla ^{h}Y_{2}\cdot Y_{3}+\nabla ^{h}\overline{r}\nabla Y_{4}\cdot
Y_{5} \\ +\left( \int\limits_{k}^{\overline{k}}\nabla ^{h}rd\kappa \right)
\nabla Y_{6}\cdot Y_{7}+\left( \int\limits_{k}^{\overline{k}}\nabla
^{h}\overline{r}d\kappa \right) \nabla ^{h}Y_{8}\cdot Y_{9} \\ +\left(
r_{z}-\int\limits_{k}^{\overline{k}}r_{z}d\kappa \right) Y_{10}+\left(
\overline{r}_{z}-\int\limits_{k}^{\overline{k}}\overline{r}_{z}d\kappa
\right) Y_{11}, \end{gathered}  \label{5.11}
\end{equation}%
where explicit expressions for functions $Y_{j}(\mathbf{x},k),j=1,...,11$
can be written in an obvious way. Also, it follows from those formulae as
well as from (\ref{5.52}) that $Y_{1},Y_{2},Y_{4},Y_{6}\in C_{1}^{h}$ and $%
Y_{3},Y_{5},Y_{7},Y_{9},Y_{10},Y_{11}\in C_{0}^{h}$. In addition,%
\begin{equation}
\left\{ 
\begin{array}{c}
\left\Vert Y_{1}\right\Vert _{C_{1}^{h}},\left\Vert Y_{2}\right\Vert
_{C_{1}^{h}},\left\Vert Y_{4}\right\Vert _{C_{1}^{h}},\left\Vert
Y_{6}\right\Vert _{C_{1}^{h}}\leq C_{1}, \\ 
\left\Vert Y_{3}\right\Vert _{C_{0}^{h}},\left\Vert Y_{5}\right\Vert
_{C_{0}^{h}},\left\Vert Y_{7}\right\Vert _{C_{0}^{h}},\left\Vert
Y_{9}\right\Vert _{C_{0}^{h}},\left\Vert Y_{10}\right\Vert
_{C_{0}^{h}},\left\Vert Y_{11}\right\Vert _{C_{0}^{h}}\leq C_{1}.%
\end{array}%
\right.   \label{5.12}
\end{equation}%
The term $Z_{1}(r,k)$ in (\ref{5.10}) is nonlinear with respect to $r$.
Applying the Cauchy-Schwarz inequality and (\ref{5.53}), we obtain%
\begin{equation}
Z_{1}(r,k)\geq \frac{1}{2}|\Delta ^{h}r|^{2}-C_{1}|\nabla
^{h}r|^{2}-C_{1}\int\limits_{k}^{\overline{k}}|\nabla ^{h}r|^{2}d\kappa .
\label{5.14}
\end{equation}%
Similarly with (\ref{5.10})--(\ref{5.14}) we obtain 
\begin{equation}
A_{2}=(\overline{z}_{1}-\overline{z}_{2})z_{2}=Z_{l,2}(r,k)+Z_{2}(r,k),
\label{5.15}
\end{equation}%
where the term $Z_{l,2}\left( r,k\right) $ is linear with respect to $r$ and
has the form similar with the one in (\ref{5.11}), although with different
functions $Y_{j},$ which still satisfy direct analogs of estimates (\ref%
{5.12}). As to the term $Z_{2}\left( r,k\right) ,$ it is nonlinear with
respect to $r$ and, as in (\ref{5.14}), 
\begin{equation}
Z_{2}(r,k)\geq -C_{1}|\nabla ^{h}r|^{2}-C_{1}\int\limits_{k}^{\overline{k}%
}|\nabla ^{h}r|^{2}d\kappa .  \label{5.16}
\end{equation}%
In addition, the following upper estimate is valid 
\begin{equation}
|Z_{1}(r,k)|+|Z_{2}(r,k)|\leq C_{1}\left( |\Delta ^{h}r|^{2}+|\nabla
^{h}r|^{2}+\int\limits_{k}^{\overline{k}}|\nabla ^{h}r|^{2}d\kappa \right) .
\label{5.17}
\end{equation}

Thus, it follows from (\ref{500.4}) and (\ref{5.7})-(\ref{5.16}) that 
\begin{equation}
\begin{gathered} J_{\lambda }(p+r)-J_{\lambda }(p)= \\ e^{2\lambda
d}\sum\limits_{j,s=1}^{N_{h}}h^{2}\int\limits_{\underline{k}}^{\overline{k}}%
\int\limits_{-\xi }^{d}\left[ (S_{1}\Delta ^{h}r+S_{2}\cdot \nabla
^{h}r)(x_{j},y_{s},z)\right] \varphi _{\lambda }(z)dzd\kappa \\ +e^{2\lambda
d}\sum\limits_{j,s=1}^{N_{h}}h^{2}\int\limits_{\underline{k}}^{\overline{k}}%
\int\limits_{-\xi }^{d}Z(r,k)(x_{j},y_{s},z)\varphi _{\lambda }(z)dzd\kappa
, \end{gathered}  \label{5.18}
\end{equation}%
where 
\begin{equation}
Z(r,k)(x_{j},y_{s},z)=Z_{1}(r,k)(x_{j},y_{s},z)+Z_{2}(r,k)(x_{j},y_{s},z).
\label{5.180}
\end{equation}%
The second line of (\ref{5.18}),%
\begin{equation}
\mathrm{Lin}(r)=e^{2\lambda d}\sum\limits_{j,s=1}^{N_{h}}h^{2}\int\limits_{%
\underline{k}}^{\overline{k}}\int\limits_{-\xi }^{d}\left[ (S_{1}\Delta
^{h}r+S_{2}\cdot \nabla ^{h}r)(x_{j},y_{s},z)\right] \varphi _{\lambda
}(z)dzd\kappa   \label{5.19}
\end{equation}%
is linear with respect to $r$. Also, the vector functions $S_{1}(\mathbf{x}%
,k)$ and $S_{2}(\mathbf{x},k)$ are such that 
\begin{equation}
\left\vert S_{1}(\mathbf{x},k)\right\vert ,\left\vert S_{2}(\mathbf{x}%
,k)\right\vert \leq C_{1}\quad \text{ in }\overline{\Omega _{h}}\times
\lbrack \underline{k},\overline{k}].  \label{5.20}
\end{equation}%
As to the third line of (\ref{5.18}), it can be estimated from the below as%
\begin{align}
& e^{2\lambda d}\sum\limits_{j,s=1}^{M_{h}}h^{2}\int\limits_{\underline{k}}^{%
\overline{k}}\int\limits_{-\xi }^{d}Z(r,k)(x_{j},y_{s},z)\varphi _{\lambda
}(z)dzd\kappa   \notag \\
& \geq e^{2\lambda d}\sum\limits_{j,s=1}^{M_{h}}h^{2}\left[ \frac{1}{2}%
\int\limits_{\underline{k}}^{\overline{k}}\int\limits_{-\xi }^{d}|\Delta
^{h}r(x_{j},y_{s},z)|^{2}\varphi _{\lambda }(z)dzd\kappa \right. 
\label{5.21} \\
& -\left. C_{1}e^{2\lambda d}\int\limits_{\underline{k}}^{\overline{k}%
}\int\limits_{-\xi }^{d}|\nabla ^{h}r(x_{j},y_{s},z)|^{2}\varphi _{\lambda
}\left( z\right) dzd\kappa \right] .  \notag
\end{align}%
In addition, using (\ref{5.17}) and (\ref{5.180}), we obtain 
\begin{equation}
\begin{gathered} e^{2\lambda
d}\sum\limits_{j,s=1}^{N_{h}}h^{2}\int\limits_{\underline{k}}^{\overline{k}}%
\int\limits_{-\xi }^{d}|Z(r,k)(x_{j},y_{s},z)|\varphi _{\lambda
}(z)dzd\kappa \\ \leq C_{1}e^{2\lambda
d}\sum\limits_{j,s=1}^{N_{h}}h^{2}\int\limits_{\underline{k}}^{\overline{k}}%
\int\limits_{-\xi }^{d}\Biggl( |\Delta ^{h}r|^{2}+|\nabla ^{h}r|^{2} \\
+\int\limits_{k}^{\overline{k}}|\nabla ^{h}r|^{2}d\kappa \Biggr)
(x_{j},y_{s},z)\varphi _{\lambda }(z) dz d\kappa . \end{gathered}
\label{5.22}
\end{equation}%
The functional $\mathrm{Lin}(h)$ in (\ref{5.19}) is linear with respect to $r
$. Also, by (\ref{5.19}) and (\ref{5.20}) 
\begin{equation*}
\left\vert \mathrm{Lin}(r)\right\vert \leq C_{1}e^{2\lambda \left( d+\xi
\right) }\left\Vert r\right\Vert _{H_{2}^{h}},\quad \text{for all }r\in
H_{0,2}^{h}.
\end{equation*}%
Hence, $\mathrm{Lin}(r):H_{0,2}^{h}\rightarrow \mathbb{R}$ is a bounded
linear functional. Hence, by Riesz theorem for each pair $\lambda >0$ there
exists a vector function $X_{\lambda }\in H_{0,2}^{h}$ independent on $r$
such that 
\begin{equation}
\mathrm{Lin}(r)=[X_{\lambda },r],\quad \text{for all }r\in H_{0,2}^{h}.
\label{5.23}
\end{equation}%
In addition, (\ref{5.18}), (\ref{5.22}) and (\ref{5.23}) imply that%
\begin{equation}
\left\vert J_{\lambda }(p_{1}+r)-J_{\lambda }(p_{1})-[X_{\lambda
},r]\right\vert \leq C_{1}e^{2\lambda \left( d+\xi \right) }\left\Vert
r\right\Vert _{H_{2}^{h}}^{2},\quad \text{for all }r\in H_{0,2}^{h}.
\label{5.24}
\end{equation}%
Thus, (\ref{5.18})--(\ref{5.24}) imply that $X_{\lambda }\in H_{0,2}^{h}$ is
the Frech\'{e}t derivative of the functional $J_{\lambda }(p)$ at the point $%
p,$ i.e. $X_{\lambda }=J_{\lambda }^{\prime }(p_{1})$.

Next, using (\ref{5.18}) and (\ref{5.21}), we obtain 
\begin{equation*}
\begin{gathered} J_{\lambda }(p+r)-J_{\lambda }(p)-J_{\lambda }^{\prime
}(p)(r) \\ \geq e^{2\lambda d}\sum\limits_{j,s=1}^{N_{h}}h^{2}\left[
\frac{1}{2}\int\limits_{\underline{k}}^{\overline{k}}\int\limits_{-\xi
}^{d}|\Delta ^{h}r(x_{j},y_{s},z)|^{2}\varphi _{\lambda }(z)dzd\kappa
-C_{1}\int\limits_{\underline{k}}^{\overline{k}}\int\limits_{-\xi
}^{d}|\nabla ^{h}r(x_{j},y_{s},z)|^{2}\varphi _{\lambda }(z)dzd\kappa
\right] . \end{gathered}
\end{equation*}%
We now apply Carleman estimate of Theorem 7.1 for $\lambda \geq \lambda _{0},
$%
\begin{equation*}
\begin{gathered} e^{2\lambda d}\sum\limits_{j,s=1}^{N_{h}}h^{2}\left[
\frac{1}{2}\int\limits_{\underline{k}}^{\overline{k}}\int\limits_{-\xi
}^{d}\left\vert \Delta ^{h}r(x_{j},y_{s},z)\right\vert ^{2}\varphi _{\lambda
}(z)dzd\kappa
-C_{1}\int\limits_{\underline{k}}^{\overline{k}}\int\limits_{-\xi
}^{d}\left\vert \nabla ^{h}r(x_{j},y_{s},z)\right\vert ^{2}\varphi _{\lambda
}(z)dzd\kappa \right] \\ \geq e^{2\lambda d}\sum\limits_{j,s=1}^{N_{h}}h^{2}
\left[ \int\limits_{\underline{k}}^{\overline{k}}\int\limits_{-\xi
}^{d}\left\vert r_{zz}(x_{j},y_{s},z)\right\vert ^{2}\varphi _{\lambda
}(z)dzd\kappa +C \lambda \int\limits_{-\xi }^{d}\left[
r_{z}(x_{j},y_{s},z)\right] ^{2}\varphi _{\lambda }(z)dz\right. \\ \left.
+\lambda ^{3} \int\limits_{-\xi }^{d} \left[ r(x_{j},y_{s},z)\right]
^{2}\varphi _{\lambda }(z)dz\right] -C_{1}e^{2\lambda
d}\int\limits_{\underline{k}}^{\overline{k}}\int\limits_{-\xi }^{d}|\nabla
^{h}r(x_{j},y_{s},z)|^{2}\varphi _{\lambda }(z)dzd\kappa . \end{gathered}
\end{equation*}%
Hence, from these two equations it follows that for sufficiently large $%
\lambda _{1}$ 
\begin{equation*}
\lambda _{1}=\lambda _{1}(\Omega _{h},R,\left\Vert F_{\ast }\right\Vert
_{H_{3}^{h}},\Vert V_{\ast }\Vert _{H^{2,h}(\Omega _{h})},\underline{k},%
\overline{k})\geq \lambda _{0}
\end{equation*}%
and for all $\lambda \geq \lambda _{1}$%
\begin{equation*}
\begin{gathered} J_{\lambda }(p_{1}+r)-J_{\lambda }(p_{1})-J_{\lambda
}^{\prime }(p_{1})(r) \geq e^{2\lambda d} \sum\limits_{j,s=1}^{N_{h}}h^{2}
\left[ \int\limits_{\underline{k}}^{\overline{k}}\int\limits_{-\xi
}^{d}\left\vert r_{zz}(x_{j},y_{s},z)\right\vert ^{2}\varphi _{\lambda
}(z)dzd\kappa \right. \\ \left. +C_{1} \lambda \int\limits_{-\xi }^{d}[
r_{z}(x_{j},y_{s},z)] ^{2}\varphi _{\lambda }(z)dz + C_1 \lambda
^{3}\int\limits_{-\xi}^{d}[ r(x_{j},y_{s},z)]^{2}\varphi _{\lambda
}(z)dz\right] \geq C_{1}\left\Vert r\right\Vert _{H_{2}^{h}}^{2},
\end{gathered}
\end{equation*}%
which establishes (\ref{6.4}). $\square $

\textbf{Theorem 7.4}. \emph{Assume that the conditions of Theorems 7.2 and
7.3 regarding the tail function }$V=V_{\mu \left( \delta \right) }$\emph{\
and\ the functions }$F$ and $F_{\ast }$\emph{\ are satisfied. Then the Frech%
\'{e}t derivative }$J_{\lambda }^{\prime }$\emph{\ of the functional }$%
J_{\lambda }$\emph{\ satisfies the Lipschitz continuity condition in any
ball }$B(R^{\prime })$\emph{\ as in (\ref{3.10}) with an arbitrary }$%
R^{\prime }>0.$\emph{\ More precisely, the following inequality holds with
the constant }$M=M(\Omega _{h},R^{\prime },\left\Vert F_{\ast }\right\Vert
_{H_{3}^{h}},\left\Vert V_{\ast }\right\Vert _{H^{2,h}(\Omega _{h})},\lambda
,\underline{k},\overline{k})>0$ \emph{\ depending only on listed parameters:}
\begin{equation*}
\left\Vert J_{\lambda }^{\prime }(p_{1})-J_{\lambda }^{\prime
}(p_{2})\right\Vert _{H_{2}^{h}}\leq M\left\Vert p_{1}-p_{2}\right\Vert
_{H_{2}^{h}},\quad \text{for all }p_{1},p_{2}\in B(R^{\prime }).
\end{equation*}

The proof of this theorem is completely similar with that of theorem 3.1 of 
\cite{BakushinskiiKlibanov17} and is, therefore, omitted.

Denote $P_{\overline{B}}:H_{0,2}^{h}\rightarrow \overline{B(R)}$ the
projection operator of the Hilbert space $H_{0,2}^{h}$ on $\overline{B(R)}%
\subset H_{0,2}^{h}.$ Let $p_{0}\in B(R)$ be an arbitrary point of the ball $%
B(R)$. Let the number $\gamma \in (0,1)$. Consider the following sequence: 
\begin{equation}
p_{n}=P_{\overline{B}}(p_{n-1}-\gamma J_{\lambda }^{\prime }(p_{n-1})),\text{
}n=1,2,\dots   \label{4.9}
\end{equation}%
The following theorem follows immediately from the combination of Theorems
7.3 and 7.4 with lemma 2.1 and Theorem 2.1 of \cite{BakushinskiiKlibanov17}.

\textbf{Theorem 7.5.} \emph{Assume that conditions of Theorems 7.2 and 7.3
are satisfied. Let }$\lambda \geq \lambda _{1},$ \emph{where }$\lambda _{1}$%
\emph{\ is defined in Theorem 7.3. Then there exists unique minimizer }$%
p_{\min ,\lambda }\in \overline{B(R)}$\emph{\ of the functional }$J_{\lambda
}(p)$\emph{\ on the set }$\overline{B(R)}$\emph{\ and } 
\begin{equation}
J_{\lambda }^{\prime }(p_{\min ,\lambda })(y-p_{\min ,\lambda })\geq 0,\quad 
\text{\emph{for all} }y\in H_{0,2}^{h}.  \label{4.6}
\end{equation}%
\emph{Also, there exists a sufficiently small number }$\gamma _{0}=\gamma
_{0}(\Omega _{h},R,\left\Vert F_{\ast }\right\Vert _{H_{3}^{h}},\left\Vert
V_{\ast }\right\Vert _{H^{2,h}(\Omega _{h})},\underline{k},\overline{k}%
,\lambda )\in (0,1)$\emph{\ depending only on listed parameters such that
for any }$\gamma \in (0,\gamma _{0})$\emph{\ the sequence (\ref{4.9})
converges }$p_{\min ,\lambda },$ 
\begin{equation}
\left\Vert p_{\min ,\lambda }-p_{n}\right\Vert _{H_{2}^{h}}\leq \theta
^{n}\left\Vert p_{\min ,\lambda }-p_{0}\right\Vert _{H_{2}^{h}},\quad
n=1,2,\dots   \label{4.90}
\end{equation}%
\emph{\ where the number }$\theta =\theta (\Omega _{h},R,\left\Vert F_{\ast
}\right\Vert _{H_{3}^{h}},\left\Vert V_{\ast }\right\Vert _{H^{2,h}(\Omega
_{h})},\underline{k},\overline{k},\lambda ,\gamma )\in (0,1)$ \emph{depends
only on listed parameters.}

Thus, (\ref{4.90}) estimates the convergence rate of the sequence (\ref{4.9}%
) to the minimizer $p_{\min ,\lambda }$. We now need to estimate the
convergence rate of this sequence to the exact solution. To do this, we
follow the Tikhonov regularization concept \cite{Bak,T} in Theorem 7.6 via
assuming that the exact solution $p_{\ast }\in B(R).$

\textbf{Theorem 7.6.} \emph{Assume that conditions of Theorems 7.2 and 7.3
are satisfied. Let }$\lambda _{1}$\emph{\ be the number of Theorem 7.3, }$%
\delta _{1}\in (0,e^{-4\left( d+\xi \right) \lambda _{1}})$ and $\delta \in
(0,\delta _{1}).$ \emph{Set }$\lambda =\lambda (\delta )=\ln (\delta
^{-1/(4\left( d+\xi \right) )})>\lambda _{1}.$ \emph{Furthermore, assume
that the function }$p_{\ast }\in B(R)$\emph{\ . Then there exists a number }%
\begin{equation*}
C_{2}=C_{2}(\Omega _{h},R,\left\Vert F_{\ast }\right\Vert
_{H_{3}^{h}},\left\Vert V_{\ast }\right\Vert _{H^{2,h}(\Omega _{h})},%
\underline{k},\overline{k})>0
\end{equation*}%
\emph{\ depending only on listed parameters such that}%
\begin{align}
\left\Vert p_{\ast }-p_{\min ,\lambda (\delta )}\right\Vert _{H_{2}^{h}}&
\leq C_{2}\delta ^{1/4},  \label{4.7} \\
\left\Vert c_{\ast }-c_{\min ,\lambda \left( \delta \right) }\right\Vert
_{L_{2}^{h}(\Omega _{h})}& \leq C_{2}\delta ^{1/4},  \label{4.8}
\end{align}%
\emph{In addition, the following convergence estimates hold}%
\begin{align}
\left\Vert p_{\ast }-p_{n}\right\Vert _{H_{2}^{h}}& \leq C_{2}\delta
^{1/4}+\theta ^{n}\left\Vert p_{\min ,\lambda (\delta )}-p_{0}\right\Vert
_{H_{2}^{h}},\text{ }n=1,2,\dots   \label{4.10} \\
\left\Vert c_{\ast }-c_{n}\right\Vert _{L_{2}^{h}(\Omega _{h})}& \leq
C_{2}\delta ^{1/4}+\theta ^{n}\left\Vert p_{\min ,\lambda \left( \delta
\right) }-p_{0}\right\Vert _{H_{2}^{h}},\text{ }n=1,2,\dots   \label{4.11}
\end{align}%
\emph{where }$\theta \in (0,1)$ \emph{is the number of Theorem 7.5 and
functions }$c_{\min ,\lambda (\delta )}(\bm x)$\emph{\ and }$c_{n}(\mathbf{x}%
)$\emph{\ is reconstructed from functions }$p_{\min ,\lambda (\delta )}$%
\emph{\ and }$p_{n}(\mathbf{x},k)$\emph{\ respectively using (\ref{3.20})--(%
\ref{3.30}) and (\ref{3.9}).}

\begin{remark}
Since $R>0$ is an arbitrary number and $p_{0}$ is an arbitrary point of the
ball $B\left( R\right) $, then Theorems 7.5 and 7.6 ensure the global
convergence of the gradient projection method for our case, see section 1.
We note that if a functional is non convex, then the convergence of a
gradient-like method of its minimization can be guaranteed only if the
starting point of iterations is located in a sufficiently small neighborhood
of its minimizer.
\end{remark}

\textit{Proof.} We temporarily denote $J_{\lambda }(p):=$ $J_{\lambda }(p+F)$%
, see (\ref{500.4}). We have 
\begin{equation}
J_{\lambda }(p_{\ast }+F_{\ast })=e^{2\lambda
d}\sum\limits_{j,s=1}^{N_{h}}h^{2}\int\limits_{\underline{k}}^{\overline{k}%
}\int\limits_{-\xi }^{d}|L^{h}\left( p_{\ast }+F_{\ast }\right)
(x_{j},y_{s},z,\kappa )|^{2}\varphi _{\lambda }^{2}(z)dzd\kappa =0.
\label{5.27}
\end{equation}%
It follows from (\ref{400.5}), (\ref{400.6}), (\ref{3.8})-(\ref{500.4}), (%
\ref{6.2}) (\ref{5.27}) that 
\begin{equation}
J_{\lambda }(p_{\ast }+F)\leq C_{2}\delta e^{2\lambda \left( d+\xi \right) }.
\label{5.28}
\end{equation}%
Next, using (\ref{6.4}) and , we obtain 
\begin{equation*}
J_{\lambda }(p_{\ast }+F)-J_{\lambda }(p_{\min ,\lambda (\delta
)}+F)-J_{\lambda }^{\prime }(p_{\min ,\lambda (\delta )}+F)(p_{\ast
}-p_{\min ,\lambda })\geq C_{2}\left\Vert p_{\ast }-p_{\min ,\lambda
}\right\Vert _{H_{2}^{h}}^{2}.
\end{equation*}%
Hence, since $-J_{\lambda }(p_{\min ,\lambda (\delta )}+F)\leq 0$ and by (%
\ref{4.6}) $-J_{\lambda }^{\prime }(p_{\min ,\lambda }+F)(p_{\ast }-p_{\min
,\lambda (\delta )})\leq 0,$ we obtain, using (\ref{5.28}) and recalling
that $\lambda =\ln (\delta ^{-1/(4\left( d+\xi \right) )})$:%
\begin{equation*}
\left\Vert p_{\ast }-p_{\min ,\lambda (\delta )}\right\Vert
_{H_{2}^{h}}^{2}\leq C_{2}\sqrt{\delta },
\end{equation*}%
which implies (\ref{4.7}). Estimate (\ref{4.8}) follows immediately from (%
\ref{3.20})-(\ref{3.30}), (\ref{3.9}), (\ref{6.2}) and (\ref{4.7}).

We now prove (\ref{4.10}) and (\ref{4.11}). Using (\ref{4.90}), (\ref{4.7})
and the triangle inequality, we obtain for $n=1,2,\dots $ 
\begin{align*}
\left\Vert p_{\ast }-p_{n}\right\Vert _{H_{2}^{h}}& \leq \left\Vert p_{\ast
}-p_{\min ,\lambda \left( \delta \right) }\right\Vert
_{H_{2}^{h}}+\left\Vert p_{\min ,\lambda \left( \delta \right)
}-p_{n}\right\Vert _{H_{2}^{h}}\leq C_{2}\delta ^{1/4}+\left\Vert p_{\min
,\lambda \left( \delta \right) }-p_{n}\right\Vert _{H_{2}^{h}} \\
& \leq C_{2}\delta ^{1/4}+\theta ^{n}\left\Vert p_{\min ,\lambda
}-p_{0}\right\Vert _{H_{2}^{h}},
\end{align*}%
which proves (\ref{4.10}). Next, using (\ref{3.20})-(\ref{3.30}), (\ref{3.9}%
), (\ref{6.2}), (\ref{4.90}) and (\ref{4.8}), we obtain 
\begin{equation*}
\begin{gathered} \left\Vert c_{\ast }-c_{n}\right\Vert _{L_{2}^{h}\left(
\Omega _{h}\right) } \leq \left\Vert c_{\ast }-c_{\min ,\lambda \left(
\delta \right) }\right\Vert _{L_{2}^{h}\left( \Omega _{h}\right)
}+\left\Vert c_{\min ,\lambda \left( \delta \right) }-c_{n}\right\Vert
_{L_{2}^{h}\left( \Omega _{h}\right) } \\ \leq C_{2}\delta
^{1/4}+C_{2}\left\Vert p_{\min ,\lambda \left( \delta \right)
}-p_{n}\right\Vert _{H_{2}^{h}}\leq C_{2}\delta ^{1/4}+C_{2}\left\Vert
p_{\min ,\lambda \left( \delta \right) }-p_{n}\right\Vert _{H_{2}^{h}} \\
\leq C_{2}\delta ^{1/4}+C_{2}\theta ^{n}\left\Vert p_{\min ,\lambda \left(
\delta \right) }-p_{0}\right\Vert _{H_{2}^{h}}. \end{gathered}
\end{equation*}%
The latter establishes (\ref{4.11}). $\square $

\section{Numerical Study}

\label{sec:8}

We present in this section a numerical study of the application of our
convexification method to microwave experimental backscatter data for buried
objects. One of possible applications is in the standoff detection of
explosives. We note that these data were treated in~\cite%
{KlibanovLiem17buried} by a different globally convergent method. We first
describe very briefly the measured data and its preprocessing which is
important for the application of our convexification method. We refer to~%
\cite{KlibanovLiem17buried} for all the details of data collection and
preprocessing.

\subsection{Measured data and its processing}

\label{sec:8.1}

\begin{figure}[h!]
\begin{center}
\includegraphics[width=0.6\textwidth]{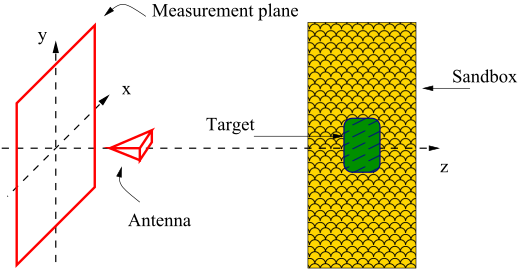}
\end{center}
\caption{A schematic diagram of the collection of our experimental data.}
\label{fig:setup}
\end{figure}

The experimental data were measured by a scattering facility at the
University of North Carolina at Charlotte. We have measured the backscatter
data for objects buried in a sandbox. This sandbox was filled with dry sand
and contains no moisture, see Figure~\ref{fig:setup}. The data were measured
on a rectangular surface of dimensions 1 m $\times $ 1 m. The distance
between this surface and the sandbox was about 75 centimeters (cm). The
coordinate system is chosen in such a way that the $x-$axis and the $y-$axis
are respectively the horizontal and the vertical axis, while the $z-$axis is
orthogonal to the measurement surface. The direction from the measurement
surface to the target is the positive direction of the $z-$axis.

The measurements consist of multi-frequency backscatter data associated with
300 frequency points uniformly distributed over the range from 1 GHz to 10
GHz. However, we work with the preprocessed data which are stable on narrow
intervals of frequencies centered at 2.6 GHz, 3.01 GHz or 3.1 GHz. Since the
corresponding wavelength for 2.6 GHz is 11.5 cm, the distance between the
source and the buried targets was about at least 6.17 wavelengths. This
distance is sufficiently large in terms of wavelengths, and therefore
justifies our modeling of the source as a plane wave. The backscatter data
were generated by a single direction of the incident plane wave.

Recall that these experimental data were preprocessed in~\cite%
{KlibanovLiem17buried} and we will study the performance of our inversion
method on that preprocessed data instead of the raw ones. The preprocessing
developed in the cited paper comprises two main goals: distill the signals
reflected by our buried targets from signals reflected by the sandbox and
other unwanted objects, and reduce the noise in the data as well as the
computational domain.

For the convenience of the readers we briefly summarize the main steps of
the data preprocessing developed in~\cite{KlibanovLiem17buried}.

\begin{enumerate}
\item[Step 1.] Subtract the reference data from the measured data for buried
objects. The reference data are the ones measured in the case when the
sandbox contains no buried objects. This subtraction helps us to sort of
extract the signals of the buried targets from the total signal and also to
reduce the noise.

\item[Step 2.] The data obtained after Step 1 were back propagated to the
sandbox using the data propagation process. This process aims to
\textquotedblleft move" the data closer to the target. As a result, we
obtain reasonable estimates for the location of the buried targets,
particularly in the $(x,y)-$plane, see Figure 1. In addition, this step
helps us reduce the computational domain.

\item[Step 3.] Determine an interval of frequencies on which the data
obtained after Step 2 are stable.
\end{enumerate}

\begin{figure}[h!]
\begin{center}
\subfloat[\label{fig:meas}]{\includegraphics[width=0.45\textwidth]{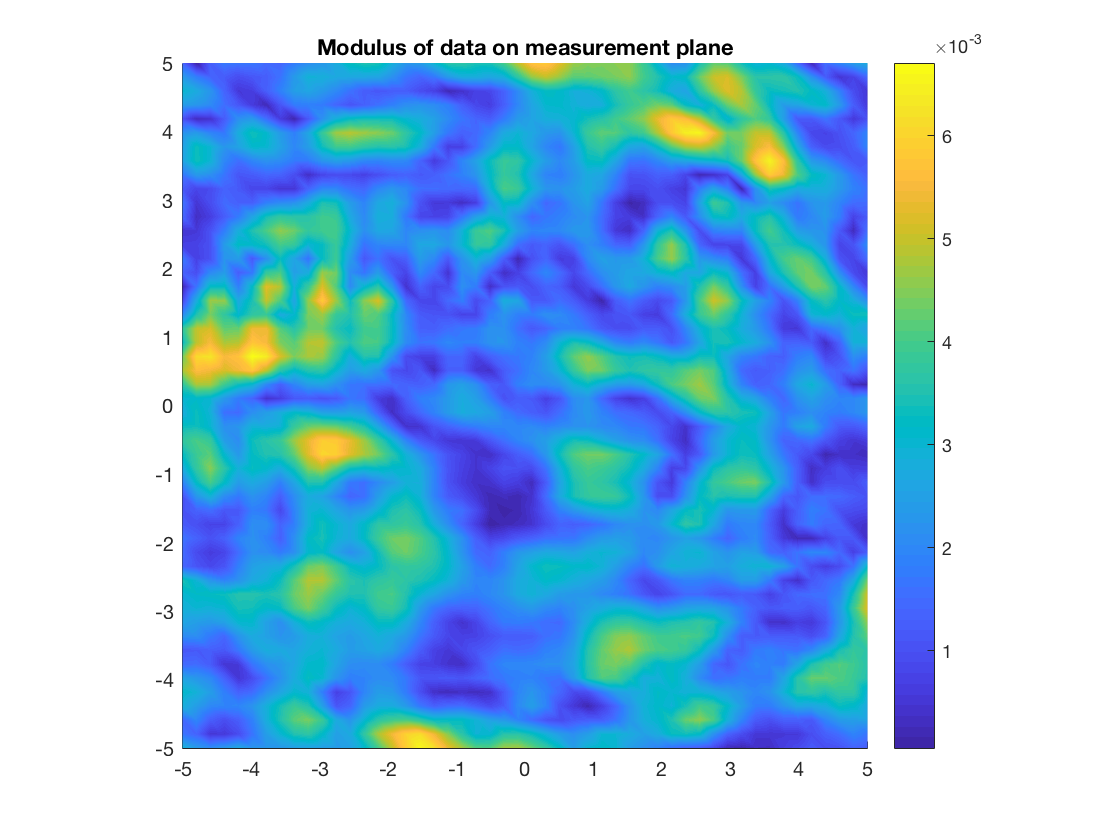}} %
\subfloat[\label{fig:prop}]{\includegraphics[width=0.45\textwidth]{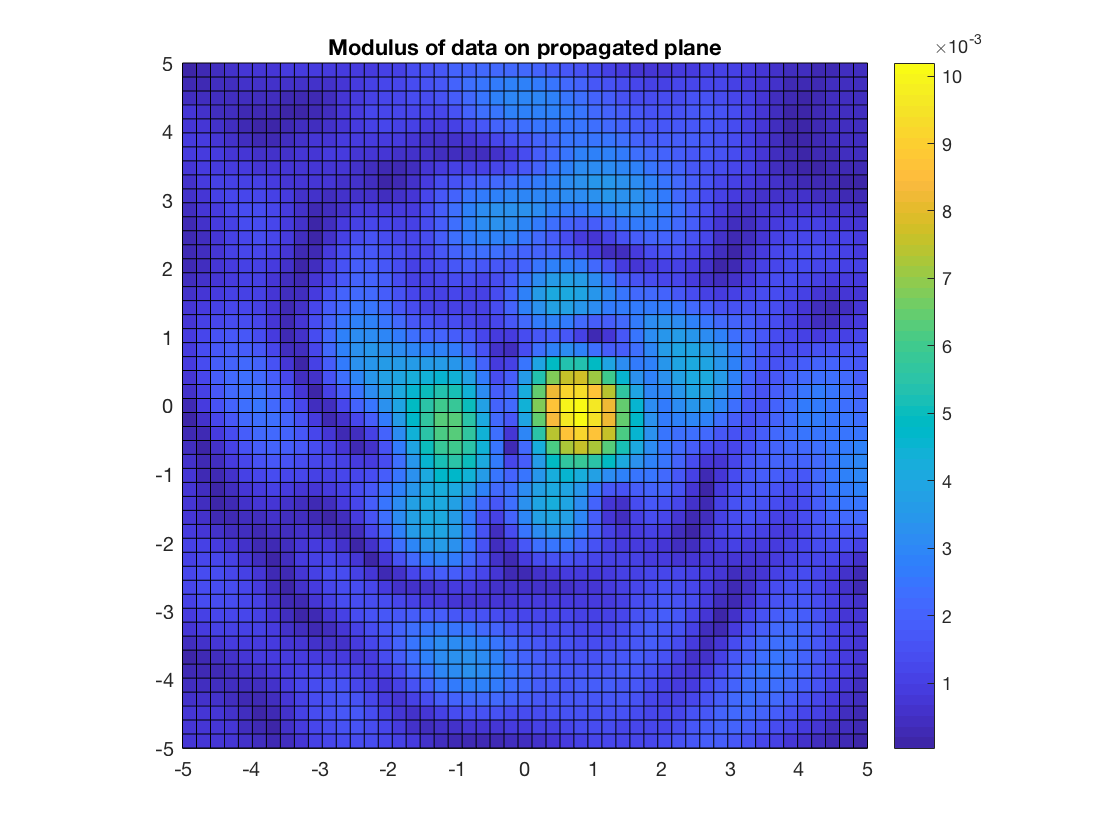}}
\end{center}
\caption{a) Absolute value of measured experimental data for the buried
target 4 (sycamore, see table \protect\ref{tab1}). b) Absolute value of the
propagated data of a).}
\label{fig:data}
\end{figure}

\subsection{Reconstruction results}

In this section we present the results of reconstructions from experimental
data for the objects buried in a sandbox in Table \ref{tab1} using our
convexification method. The experimental setup for the case of objects
buried in a sandbox is shown in Figure \ref{fig:setup}.

\begin{table}[h!]
\caption{Buried objects}
\label{tab1}
\begin{center}
\begin{tabular}{ccc}
\hline
Number & $\quad$ Description $\quad$ & $\quad$ Size in $x\times y\times z $
directions (in cm) $\quad$ \\ \hline
1 & Bamboo & $3.8 \times 11.6 \times 3.8 $ \\ 
2 & Geode & $8.8 \times 8.8 \times 8.8 $ \\ 
3 & Rock & $10.5 \times 7.5 \times 4.0 $ \\ 
4 & Sycamore & $3.8 \times 9.9 \times 3.8 $ \\ 
5 & Wet wood & $9.1 \times 5.7 \times 5.8 $ \\ 
6 & Yellow pine & $9.0 \times 8.3 \times 5.8 $ \\ \hline
\end{tabular}%
\end{center}
\end{table}

In Table \ref{tab2} we present the optimal frequencies and corresponding
intervals of wavenumbers $[\underline{k}, \overline{k}]$ for our objects. We
refer to~\cite{KlibanovLiem17buried} for the details of the determination of
these intervals.

The objects with their directly measured dielectric constant $c_{meas}$ and
computed coefficient $c_{comp}$ along with corresponding measurement $%
\varepsilon _{meas}$ and computational errors $\varepsilon
_{comp}=|c_{comp}-c_{meas}|/c_{meas}\ast 100\%$ are listed in Table \ref%
{tab3}. Note that the coefficients $c_{comp}$ in Table \ref{tab3} are the
maximal values of the reconstructed functions $c(\mathbf{x})$. In all our
numerical tests we have used reasonable values of parameters $\mu = \lambda
= 3.0.$

Considering the significant amount of noise in the measured data, the
computational errors $\varepsilon _{comp}$ of reconstructed coefficients are
sufficiently small. The computed dielectric constant of object 3 (a piece of
rock) has the biggest error $\varepsilon _{comp}=9.63\%$, but it is lower
than its measurement error $21.3\%$.

In Table \ref{tab4} we present the propagation distance $d$ \cite%
{KlibanovLiem17buried}, estimated location of objects and location of the
reconstructed objects, i.e. the location of the maximum value of computed
coefficient $\max (c_{comp}(\mathbf{x}))$. Errors of locations are small
comparable with the size of the computational domain where we solve our
inverse problem

Fig. \ref{fig:c_obj2} and \ref{fig:c_obj4} illustrate the exact and computed
images for the objects 2 and 4, respectively. Images are obtained using the
contour filter in Paraview.

\subsection{Conclusion}

\begin{figure}[tbp]
\begin{center}
\subfloat[\label{fig:c_true2}]{\includegraphics[width=0.45%
\textwidth]{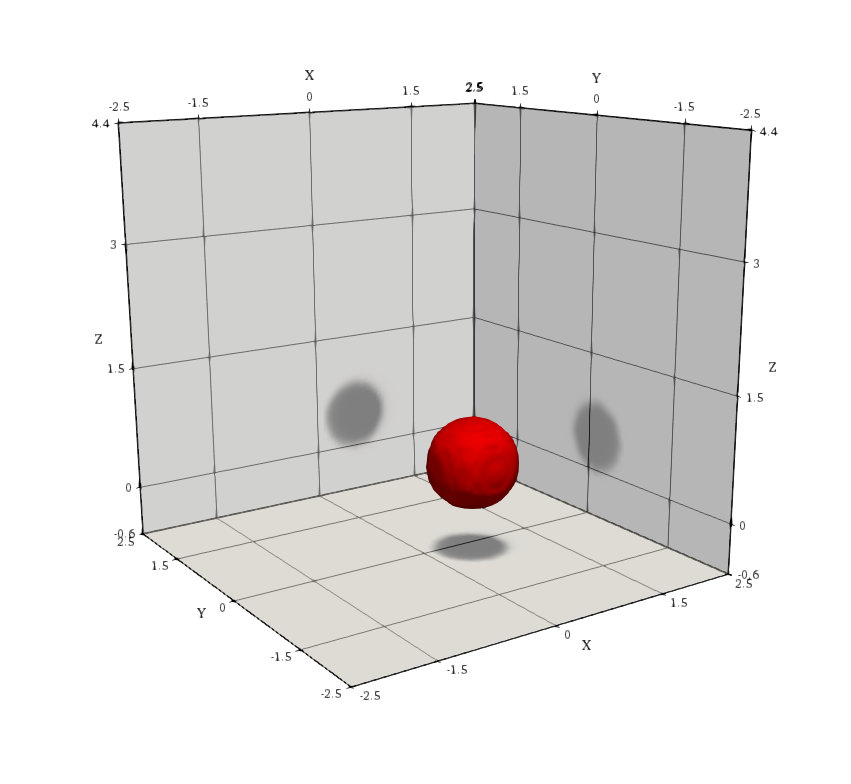}} 
\subfloat[\label{fig:c_comp2}
]{\includegraphics[width=0.45\textwidth]{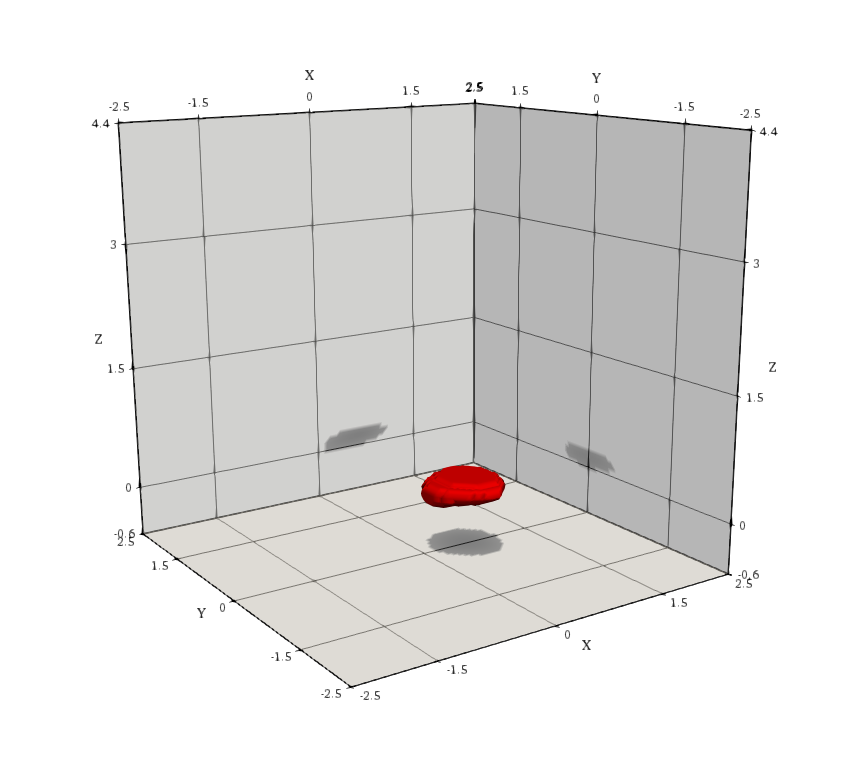}}
\end{center}
\caption{Reconstruction result for target 2: (a) exact image, (b) computed
image}
\label{fig:c_obj2}
\end{figure}

\begin{figure}[tbp]
\begin{center}
\subfloat[\label{fig:c_true4}]{\includegraphics[width=0.45%
\textwidth]{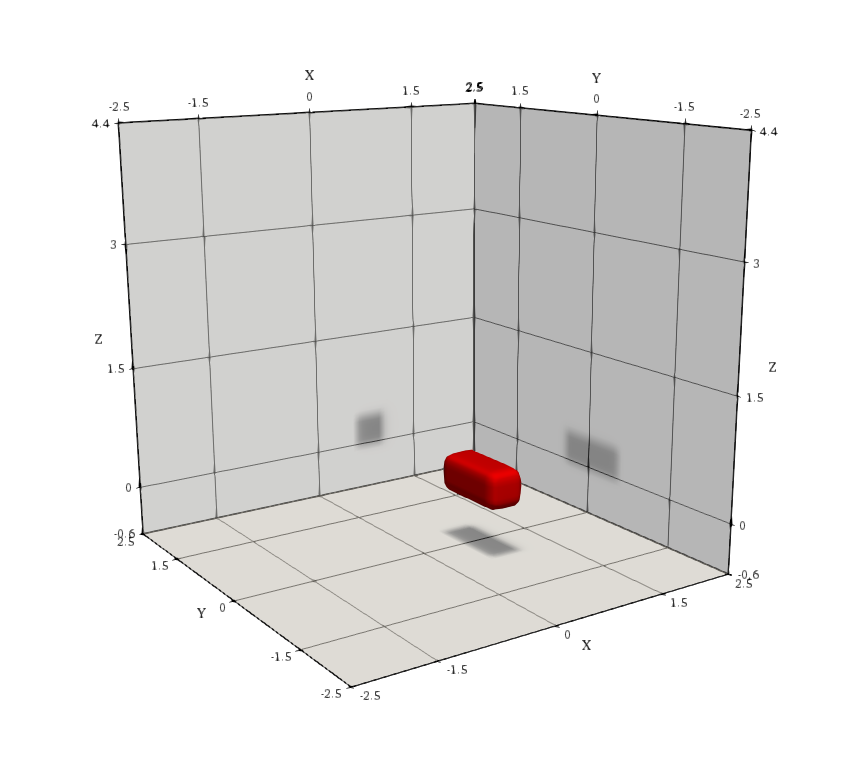}} \subfloat[\label{fig:c_comp4}]{%
\includegraphics[width=0.45\textwidth]{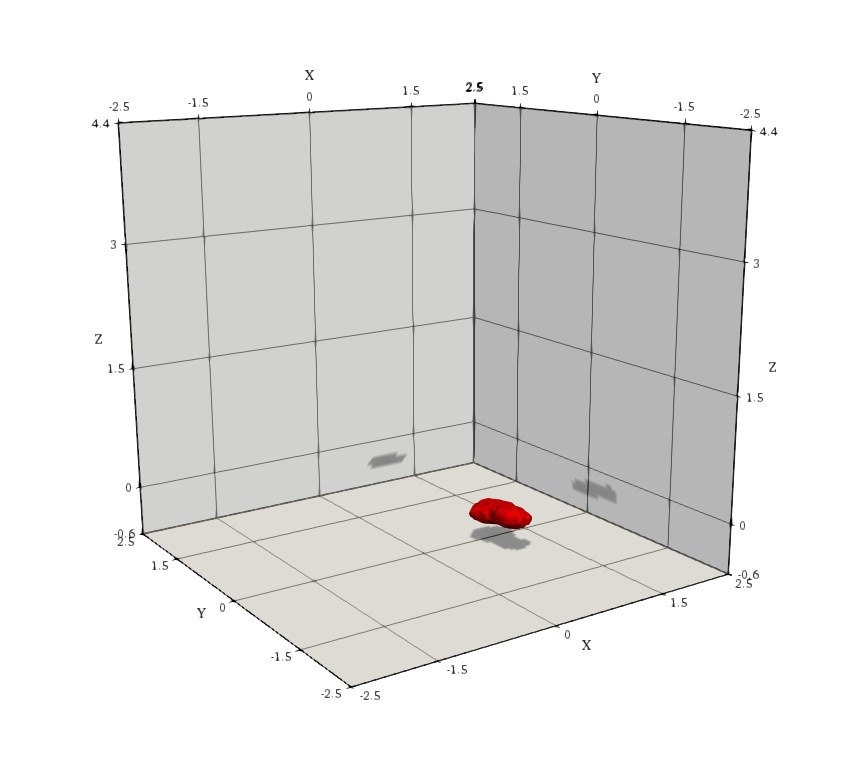}}
\end{center}
\caption{Reconstruction result for the target 4: (a) exact image, (b)
computed image}
\label{fig:c_obj4}
\end{figure}

Table \ref{tab3} and Figures \ref{fig:c_obj2}, \ref{fig:c_obj4} demonstrate
that our numerical method accurately reconstructs both dielectric constants
and locations of targets in a quite challenging case of backscatter
experimental data collected for buried targets.


\begin{table}[h!]
\caption{Optimal frequencies and interval of wavenumbers}
\label{tab2}
\begin{center}
\begin{tabular}{ccc}
\hline
Number & Optimal frequency, GHz & Interval of wavenumbers $[\underline{k},
\, \overline{k}]$ \\ \hline
1 & 3.10 & [6.322, 6.638] \\ 
2 & 3.01 & [6.133, 6.448] \\ 
3 & 3.01 & [6.070, 6.385] \\ 
4 & 3.10 & [6.322, 6.638] \\ 
5 & 2.62 & [5.313, 5.691] \\ 
6 & 2.62 & [5.313, 5.691] \\ \hline
\end{tabular}%
\end{center}
\end{table}

\begin{table}[h!]
\caption{Measured and reconstructed coefficients of objects}
\label{tab3}
\begin{center}
\begin{tabular}{ccccc}
\hline
Number & $c_{meas}$ & $\varepsilon_{meas}$ & $c_{comp}$ & $\varepsilon_{comp}
$ \\ \hline
1 & 4.50 & $5.99 \%$ & 4.69 & $4.22 \%$ \\ 
2 & 5.45 & $1.13 \%$ & 5.28 & $3.12 \%$ \\ 
3 & 5.61 & $21.3 \%$ & 5.07 & $9.63 \%$ \\ 
4 & 4.89 & $2.89 \%$ & 4.95 & $1.23 \%$ \\ 
5 & 7.58 & $4.69 \%$ & 8.06 & $6.33 \%$ \\ 
6 & 4.89 & $1.54 \%$ & 5.22 & $8.75 \%$ \\ \hline
\end{tabular}%
\end{center}
\end{table}

\begin{table}[h!]
\caption{Estimated and reconstructed locations of objects}
\label{tab4}
\begin{center}
\begin{tabular}{cccc}
\hline
Number & Estimated location in $(x,y,z)$ & Computed location in $(x,y,z)$ & 
\\ \hline
1 & (0.80, -0.11, 0.19) & (0.83, 0.03, -0.05) &  \\ 
2 & (0.58, -0.14, 0.44) & (0.63, 0.03, 0.16) &  \\ 
3 & (0.62, -0.14, 0.20) & (0.63, 0.08, -0.20) &  \\ 
4 & (0.80, -0.04, 0.19) & (1.04, 0.08, -0.30) &  \\ 
5 & (0.57, -0.42, 0.29) & (0.53, -0.08, 0.16) &  \\ 
6 & (0.54, -0.33, 0.29) & (0.53, -0.03, 0.21) &  \\ \hline
\end{tabular}%
\end{center}
\end{table}

\begin{thebibliography}{99}
\bibitem{Am1} \textsc{H. Ammari, J. Garnier, W. Jing, H. Kang, M. Lim, K.
Solna, and H. Wang}, \emph{Mathematical and Statistical Methods for
Multistatic Imaging}, vol. 2098 of Lecture Notes in Mathematics, Springer,
Cham, 2013.

\bibitem{Am2} \textsc{H. Ammari, Y.T. Chow, and J. Zou}, \emph{The concept
of heterogeneous scattering and its applications in inverse medium scattering%
}, SIAM J. Math. Anal., 46 (2014), 2905-2935.

\bibitem{BakushinskiiKlibanov17} \textsc{A.~B. Bakushinskii, M.~V. Klibanov,
and N.~A. Koshev}, \emph{Carleman weight functions for a globally convergent
numerical method for ill-posed Cauchy problems for some quasilinear PDEs},
Nonlinear Analysis: Real World Applications, 34 (2017), pp.~201--224.

\bibitem{Bak} \textsc{A.~B. Bakushinskii, M.Yu. Kokurin and M.M. Kokurin, }%
Regularization Algorithms for Ill-Posed Problems, De Guyter, Berlin, 2018.

\bibitem{Baud} \textsc{L. Baudouin, M. de Buhan and S. Ervedoza}, \emph{%
Convergent algorithm based on Carleman estimates for the recovert of a
potential in the wave equation}, SIAM J. on Numerical Analysis, 55 (2017),
1578-1613.

\bibitem{BeilinaKlibanov12} \textsc{L.~Beilina and M.~V. Klibanov}, \emph{%
Approximate global convergence and adaptivity for coefficient inverse
problems}, Springer, 2012.


\bibitem{BY} \textsc{M.~Bellassoued and M.~Yamamoto}, \emph{Carleman
Estimates and Applications to Inverse Problems for Hyperbolic Systems},
Springer Japan KK, 2017.

\bibitem{Burman} \textsc{E. Burman, J. Ish-Horowicz and L. Oksanen}, Fully
discrete finite element data assimilation \ method for the heat equation, 
\emph{arxiv}: 1707.06908, 2017.

\bibitem{KlibanovBukhgeim81} \textsc{A.~Bukhgeim and M.~Klibanov}, \emph{%
Uniqueness in the large of a class of multidimensional inverse problems},
Soviet Math. Doklady, 17 (1981), pp.~244--247.

\bibitem{Cakon2006} \textsc{F.~Cakoni and D.~Colton}, \emph{Qualitative
Methods in Inverse Scattering Theory. An Introduction}, Springer, Berlin,
2006.

\bibitem{Chavent09} \textsc{G.~Chavent}, \emph{Nonlinear Least Squares for
Inverse Problems - Theoretical Foundations and Step-by-Step Guide for
Applications}, Springer, 2009.

\bibitem{Colto2013} \textsc{D.~Colton and R.~Kress}, \emph{Inverse Acoustic
and Electromagnetic Scattering Theory}, Springer, New York, 3rd~ed., 2013.

\bibitem{Goncharsky13} \textsc{A.~Goncharsky and S.~Romanov}, \emph{%
Supercomputer technologies in inverse problems of ultrasound tomography},
Inverse Problems, 29 (2013), p.~075004.

\bibitem{Goncharsky17} \textsc{A.~V. Goncharsky and S.~Y. Romanov}, \emph{%
Iterative methods for solving coefficient inverse problems of wave
tomography in models with attenuation}, Inverse Problems, 33 (2017),
p.~025003.

\bibitem{It} \textsc{K. Ito, B. Jin, and J. Zou}, \emph{A direct sampling
method to an inverse medium scattering problem}, Inverse Problems 28 (2012),
025003.

\bibitem{Ito} \textsc{K. Ito, B. Jin, and J. Zou}, \emph{A direct sampling
method for inverse electromagnetic medium scattering}, Inverse Problems 29
(2013), 095018.

\bibitem{Kab1} \textsc{S. I. Kabanikhin, A. D. Satybaev, M. Shishlenin}, \emph{%
Direct Methods of Solving Multidimensional Inverse Hyperbolic Problem}, VSP,
Utrecht, 2004.

\bibitem{Kab2} \textsc{S. Kabanikhin, K. Sabelfeld, N. Novikov, M. Shishlenin%
}, \emph{Numerical solution of the multidimensional Gelfand-Levitan equation}%
, J. Inverse and Ill-Posed Problems 23 (2015) 439-450.

\bibitem{KS} \textsc{M. V. Klibanov and F. Santosa}, \emph{A computational
quasi-reversibility method for Cauchy problems for Laplace's equation}, SIAM
J. Applied Mathematics, 51 (1991), pp. 1653-1675.

\bibitem{Klibanov97a} \textsc{M.~V. Klibanov}, \emph{Global convexity in a
three-dimensional inverse acoustic problem}, SIAM Journal on Mathematical
Analysis, 28 (1997), pp.~1371--1388.


\bibitem{KlibanovTimonov04} \textsc{M.~V. Klibanov and A.~Timonov}, \emph{%
Carleman Estimates for Coefficient Inverse Problems and Numerical
Applications}, de Gruyter, Utrecht, 2004.

\bibitem{Ksurvey} \textsc{M.~V. Klibanov}, \emph{Carleman estimates for
global uniqueness, stability and numerical methods for coefficient inverse
problems}, Journal of Inverse and Ill-Posed Problems, 21 (2013),
pp.~477--560.

\bibitem{KQR} \textsc{M.~V. Klibanov}, \emph{Carleman estimates for the
regularization of ill-posed Cauchy problems, }Applied Numerical Mathematics,
94 (2015), pp. 46--74.

\bibitem{Klibanov15} \textsc{M.~V. Klibanov}, \emph{Carleman weight
functions for solving ill-posed Cauchy problems for quasilinear PDEs},
Inverse Problems, 31 (2015), p.~125007.

\bibitem{KlibanovThanh15} \textsc{M.~V. Klibanov and N.~T. Th{\`{a}}nh}, 
\emph{Recovering dielectric constants of explosives via a globally strictly
convex cost functional}, SIAM Journal on Applied Mathematics, 75 (2015),
pp.~518--537.

\bibitem{KlibanovRomanov16} \textsc{M.~V. Klibanov and V.~Romanov}, \emph{%
Two reconstruction procedures for a 3-D phaseless inverse scattering problem
for the generalized Helmholtz equation}, Inverse Problems, 32 (2016),
p.~0150058.


\bibitem{KlibanovKolesov17} \textsc{M. V. Klibanov, A. E. Kolesov,
L. Nguyen, and A. Sullivan}, \emph{Globally strictly convex cost functional
for a 1-D inverse medium scattering problem with experimental data}, SIAM J.
Appl. Math., 77 (2017), 1733-1755.

\bibitem{KlibanovLiem16} \textsc{M.~V. Klibanov, D.-L. Nguyen, L.~H. Nguyen,
and H.~Liu}, \emph{A globally convergent numerical method for a 3D
coefficient inverse problem with a single measurement of multi-frequency data%
}, Inverse Problems and Imaging\emph{, }12 (2018), 493-523.

\bibitem{KlibKol3D} \textsc{M. V. Klibanov and A. E. Kolesov, }\emph{%
Convexification of a 3-D coefficient inverse scattering problem}, Computers
and Mathematics with Applications, published online,
https://doi.org/10.1016/j.camwa.2018.03.016, 2018.



\bibitem{Lakhal1} \textsc{A. Lakhal}, \emph{A decoupling-based imaging
method for inverse medium scattering for Maxwell's equations}, Inverse
Problems, 26 (2010), 015007.

\bibitem{Liu1} \textsc{J. Li, H. Liu, and Q. Wang}, \emph{Enhanced
multilevel linear sampling methods for inverse scattering problems}, J.
Comput. Phys., 257 (2014), pp. 554--571.

\bibitem{Liu2} \textsc{J. Li, P. Li, H. Liu, and X. Liu}, \emph{Recovering
multiscale buried anomalies in a two-layered medium}, Inverse Problems, 31
(2015), 105006.

\bibitem{Kar1} \textsc{L. A. Nazarova, L. A. Nazarov, A. L. Karchevsky, M.
Vandamme}, \emph{Determining kinetic parameters of a block coal bed gas by
solving inverse problem based on data of borehole gas measurements}, Journal
of Mining Science, 2015, Vol. 51, No. 4, pp. 666--672.

\bibitem{Kar2} \textsc{A. A. Duchkov, A. L. Karchevskii, }Application of
temperature monitoring to estimate the heat flux and thermophysical
properties of bottom sediments, \emph{Doklady Earth Sciences}, 2014, Vol.
458, Part 2, p. 1285--1288.


\bibitem{KlibanovLiem17buried} \textsc{D.-L. Nguyen, M.~V. Klibanov, L.~H.
Nguyen, and M.~A. Fiddy}, \emph{\ {Imaging of buried objects from
multi-frequency experimental data using a globally convergent inversion
method}}, J. Inverse and Ill-Posed Problems, accepted for publication
(2017), available online of this journal, DOI: 10.1515/jiip-2017-0047; also
available at arxiv: 1705.01219, 2017.

\bibitem{KlibanovLiem17exp} \textsc{D.-L. Nguyen, M.~V. Klibanov, L.~H.
Nguyen, A.~E. Kolesov, M.~A. Fiddy, and H.~Liu}, \emph{Numerical solution of
a coefficient inverse problem with multi-frequency experimental raw data by
a globally convergent algorithm}, Journal of Computational Physics, 345
(2017), pp.~17--32.

\bibitem{Rom3} \textsc{R. G. Romanov}, \emph{Inverse Problems of
Mathematical Physics}, VSP, Utrecht, 1986.

\bibitem{Rom} \textsc{V. G. Romanov}, \emph{Inverse problems for differential
equations with memory}, Eurasian J. of Mathematical and Computer
Applications, 2, issue 4, pp. 51-80, 2014.

\bibitem{Scales92} \textsc{J.~A. Scales, M.~L. Smith, and T.~L. Fischer}, 
\emph{Global optimization methods for multimodal inverse problems}, Journal
of Computational Physics, 103 (1992), pp.~258--268.


\bibitem{T} A.N. Tikhonov, A.V. Goncharsky, V.V. Stepanov and A.G. Yagola, 
\emph{Numerical Methods for the Solution of Ill-Posed Problems}, Kluwer,
London, 1995.

\end{thebibliography}
\end{document}